\documentclass{report}

\usepackage{amsfonts,amsmath,amssymb}
\usepackage[dvips]{graphics}
\usepackage{epsfig}
\usepackage{amscd}
\usepackage{verbatim}
\usepackage{euscript}
\usepackage{color}
\usepackage{subfigure}
\usepackage{thumbpdf}
\usepackage{hyperref}

\newtheorem{theorem}{\sc Theorem}[section]
\newtheorem{lemma}{\sc Lemma}[section]
\newtheorem{proposition}{\sc Proposition}[section]
\newtheorem{remark}{\sc Remark}[section]
\newtheorem{definition}{\sc Definition}[section]

\def\ni{\noindent}

\def\a{\alpha}

\def\lap{\Delta}

\newcommand{\rem}[1]{{\bf Remark:}}

\def\qed{{\hspace*{\fill}{\vrule height 1ex width 1ex }\quad} 
    \vskip 0pt plus20pt}
\begin{document}
\pagenumbering{roman}
\begin{center}
{\Large \bf Well-posedness of the three-dimensional Lagrangian averaged 
Navier-Stokes equations}
\vskip 0.2in
\centerline{By}
\vskip 0.1in
\centerline{\large  James P. Peirce}
\centerline{\large B.S. (University of Washington) 1997}

\vskip 0.3in
\centerline{\large  DISSERTATION}
\vskip 0.1in
\centerline{\large Submitted in partial satisfaction of the requirements for 
the degree of}
\vskip 0.1in
\centerline{\large  DOCTOR OF PHILOSOPHY}
\vskip 0.1in
\centerline{\large in}
\vskip 0.1in
\centerline{\large MATHEMATICS}
\vskip 0.2in
\centerline{\large in the}
\vskip 0.1in
\centerline{\large  OFFICE OF GRADUATE STUDIES}
\vskip 0.1in
\centerline{\large of the}
\vskip 0.1in
\centerline{\large  UNIVERSITY OF CALIFORNIA}
\vskip 0.1in
\centerline{\large  DAVIS}
\end{center}
\vskip 0.1in
{\large Approved:}
\begin{center}
\centerline{STEVE SHKOLLER{\hskip 2.5in}}
\vskip 0.15in
\centerline{ALLEN EDELSON{\hskip 2.5in}}
\vskip 0.15in
\centerline{ALBERT FANNJIANG{\hskip 2.5in}}
\vskip 0.2in
\centerline{\large Committee in Charge}
\vskip 0.2in
\centerline{\large 2004}
\end{center}

\newpage
\large
\tableofcontents


\newpage
\large
{\Large \bf ACKNOWLEDGMENTS} \\

I have been very lucky to be a part of the U.C. Davis Mathematics Department 
for the past seven years.  My advisor, Professor Steve Shkoller, has been
amazingly patient and I thank him for his guidance.  In addition, 
Professor Allan Edelson is a good friend and I greatly appreciate his 
support early in my career at U.C. Davis.  Most importantly I would like to 
thank my family, especially my father James and mother Suzanne, for their 
unconditional support.  Without them, I may have never have survived 
this experience.  I am thankful to the department and its members for 
providing a good environment for growth and learning.  In particular I 
would like to thank Professor Bruno Nachtergaele, Professor Blake Temple, 
Professor John Hunter, Dr. Wolfgang Spitzer, and Dr. Daniel Coutand for 
their help, guidance, and encouragement.


\newpage 
\begin{center}
\underline{\bf \Large Abstract} 
\end{center}

\medskip

\large
In this dissertation we study the well-posedness of the 
three-dimensional Lagrangian averaged Navier-Stokes (LANS-$\alpha$) 
equations.  The LANS-$\alpha$ equations are a system of PDEs designed to 
capture the large scale dynamics 
of the incompressible Navier-Stokes equations.  In the Lagrangian 
averaging approach, the motion at spatial scales smaller than a 
chosen parameter $\alpha>0$ is filtered without the use of 
artificial viscosity.  There are two types of 
LANS-$\alpha$ equations: the anisotropic version in which the fluctuation 
tensor is a dynamic variable that is coupled with the evolution equations 
for the mean velocity, and the isotropic version in which the covariance 
tensor is assumed to be a constant multiple of the identity matrix.  

We prove the global-in-time existence and uniqueness of weak solutions 
to the isotropic LANS-$\alpha$ equations for the case of no-slip boundary 
conditions, generalizing the known periodic box result \cite{fht02}.  Our 
proof makes use of a formulation of the equations on bounded domains 
provided by Marsden and Shkoller \cite{ms01}.  In the anisotropic 
model, there are two choices for the divergence-free projection 
of the viscosity term.  One choice is the classic Leray 
projector.  In this case, Marsden and Shkoller \cite{ms03} have shown the 
local-in-time well-posedness of the anisotropic equations in the periodic 
box.  We extend their result by considering 
the second choice of projector, the generalized Stokes projector.  The 
local-in-time well-posedness of the anisotropic LANS-$\alpha$ with 
this viscosity term is proven by using quasi-linear PDE-type methods.

We numerically compute strong solutions to the anisotropic equations in 
the laminar channel and pipe by considering steady fluid flow with no-slip 
boundary conditions.  In particular, given a steady 
velocity vector that solves the Navier-Stokes equations, we numerically 
calculate the covariance tensor such that the pair solves the anisotropic 
LANS-$\alpha$ equations.  Our solutions are in good agreement with the
results contained in \cite{cs03a}.  Namely, we confirm the logarithmic 
degeneracy rate of the covariance tensor near the boundary and show that in 
elementary domains the sup-norm of the covariance tensor is unbounded 
in time near the wall.  We conclude the dissertation by showing the existence 
of shear flow solutions to the anisotropic LANS-$\alpha$ equations.  


\newpage
\pagestyle{myheadings} 
\pagenumbering{arabic}
\markright{  \rm \normalsize CHAPTER 1. \hspace{0.5cm}
 Introduction }
\large 
\chapter{Introduction}
\label{chapter:introduction}
\thispagestyle{myheadings}
\section{The Navier-Stokes equations.}

The Navier-Stokes equations for an incompressible fluid are 
a system of partial differential equations that model the 
velocity vector $u=u(t,x)$ and pressure function $p=p(t,x)$ of a fluid whose 
velocity is divergence-free.  We 
assume the fluid is contained in a fixed domain $\Omega \subset 
\mathbb{R}^{n}$ with boundary $\delta \Omega$ (possibly empty) and has 
a constant density $\rho$ with value $\rho=1$.  In this case, the 
Navier-Stokes equations are given as
\begin{equation}\label{eq:ns}
\begin{split} 
    \frac{\partial u}{\partial t} + (u \cdot \nabla) u = &- 
        \operatorname{grad} p + \nu \lap u + f, \\
        \operatorname{div} u = 0, \,&\, u(0,x)=u_{0}(x),\\
        u=0 & \mbox{ on } \partial \Omega,
\end{split}
\end{equation}
where $\nu>0$ is the kinematic viscosity and the vector $f$ represents 
external forces.  The kinematic viscosity is defined as the 
coefficient of viscosity divided by the density of the fluid.  
The pressure function is determined (modulo a 
constant) from the incompressibility constraint $\operatorname{div} u 
=0$ by solving the Neumann problem $-\lap p = \operatorname{div} (u 
\cdot \nabla) u$ with boundary condition $\operatorname{grad} p \cdot 
n = 0$, where $n$ denotes the outward normal vector to the boundary.  In the 
limiting case of viscosity $\nu \rightarrow 0$, the Navier-Stokes 
equations are reduce to the idealized setting of the Euler equations.

The ratio of forces from the convective nonlinearity $(u \cdot \nabla)u$ 
term and the linear diffusion $\nu \lap u$ term is called the 
Reynolds number.  As the Reynolds number increases, the contribution 
to the motion of the fluid from the viscous term decreases and turbulence is 
introduced.  In 
turbulent regimes, the nonlinear effects send energy from the large spatial 
scales to smaller and smaller scales until the energy reaches the Kolmogorov 
dissipation scale, at which it is abolished by the linear dispersive 
mechanism.  To resolve a numerical simulation of the 
Navier-Stokes equations (\ref{eq:ns}), enough grid points or Fourier 
modes must be used so that the approximation captures the energy 
cascade in all scales down to the Kolmogorov scale.  For turbulent 
flows such resolution requirements are not yet achievable, making the 
problem of turbulence an important unsolved problem in 
physics.  The numerical inability in resolving small spatial scales motivates 
the study of the averaged motion of an incompressible fluid.  The averaged 
fluid methodology, discussed in the next section, captures the dynamics of the 
large scale motion while averaging the computationally unresolvable 
scales of the Navier-Stokes equations.  

\section{Averaged fluid motion.}

An approach to modeling the averaged motion of an 
incompressible fluid is to suppose that the 
velocity of the fluid is a random variable represented by 
the decomposition 
\begin{equation}\label{eq:rd}
u(t,x) = U(t,x) + u'(t,x),
\end{equation}
where $U$ denotes the mean velocity field and $u'$ is a random 
variable with mean zero.  In a statistical theory for 
turbulence, the evolution of the fluid at large spatial scales is the 
primary focus.  The process of substituting equation 
(\ref{eq:rd}) into the Navier-Stokes equations and averaging results in the 
Reynolds averaged Navier-Stokes (RANS) equations \cite{h75,t67}.  The RANS 
equations are written as 
\begin{equation}\label{eq:RANS}
    \begin{split}
\frac{\partial U}{\partial t} + (U \cdot \nabla) U + 
\operatorname{div} &\overline{( u' \otimes u' )} =  -\operatorname{grad} p + 
        \nu \lap U + f \\
        \operatorname{div} U = 0, \,&\, U(0,x)=u_{0}(x) \\
        U = 0 & \mbox{ on } \partial \Omega.
    \end{split}
\end{equation}
The dynamics of the mean velocity $U$ is well-defined when the 
Reynolds stress tensor $\overline{( u' \otimes u' )}$ is expressed 
in terms of $U$.  The effect the motion at small length scales 
has on the evolution of the averaged velocity $U$ is called the turbulence 
closure problem.  Classically, it is assumed that the 
Reynolds stress term is of the form $\displaystyle{ 
\operatorname{div}\overline{(u' \otimes u' )}} = 
\nu_{E}(t,x,\operatorname{Def} U) \cdot \operatorname{Def} U,$ where 
$\nu_{E}$ is the eddy viscosity and $\operatorname{Def} U$ is the rate of 
deformation tensor defined as
\begin{equation}\label{eq:Def}
    \operatorname{Def} U = \frac{1}{2} \left[ \nabla U + \left(\nabla 
    U \right)^{T} \right].
\end{equation}
Under this assumption, viscosity, that is not naturally present in the 
physical model, is added into the system.  This artificial viscosity 
augments the inherent dissipative mechanism and assists in the removal 
of the energy contained in the small scales at which $u'$ resides.  Since it is 
still necessary to guess the form of $\nu_{E}$, an improvement to the 
procedure of modeling the averaged motion of a fluid is needed.  

The Lagrangian averaged Navier-Stokes (LANS-$\alpha$) equations are a 
system of partial differential equations designed to capture the 
large scale dynamics of the incompressible Navier-Stokes equations 
without the use of artificial viscosity or dissipation.  In the 
Lagrangian averaging approach, the motion at spatial scales smaller than a 
chosen parameter $\alpha > 0$ are filtered.  The inviscid form of the 
LANS-$\alpha$ equations, the Lagrangian averaged Euler (LAE-$\alpha$) 
equations, first appeared in Holm, Marsden, and 
Ratiu \cite{hmr98a, hmr98b} as a $n$-dimensional generalization of the 
one-dimensional Camassa-Holm equation.  The authors expressed the LAE-$\alpha$ 
equations as the Euler-Poincar\'{e} equations corresponding to a Lagrangian 
given by an $H^{1}$-equivalent norm.  In the next section, we discuss their 
results in greater detail.  

Marsden and Shkoller provide a 
complete derivation of the LANS-$\alpha$ equation in 
\cite{ms01,ms03}.  Rather than 
averaging at the level of the Navier-Stokes equations, the authors 
average at the level of the action functional.  For the LAE-$\alpha$ equations, 
the action functional is an $\alpha^{2}$-modification to the kinetic 
energy action functional of the classical Euler equations.  Unlike the
Reynolds averaging procedure described above, averaging at the 
level of the action functional preserves the variational structure of 
the incompressible fluid.  In particular, the solutions to the LAE-$\alpha$ 
equation are the extrema of an $H^{1}$-equivalent energy functional just as 
the solutions to the Euler equations are minimizers of the total kinetic energy.   
To highlight additional features of the 
Lagrangian averaging approach, we outline Marsden and Shkoller's recent 
derivation \cite{ms03}.  For a chosen parameter $\alpha > 0$, the authors  
define the initial data $u^{\varepsilon}_{0} = u_{0}+\varepsilon 
\omega$, $\omega \in S^{2}$, in a ball of radius $0<\varepsilon<\alpha$ 
centered at given initial data $u_{0}$.  Solving the Euler 
equations with initial data 
$u_{0}$ and $u_{0}^{\varepsilon}$ results in a velocity field $u$ and 
a perturbed velocity field $u^{\varepsilon}$.  Define $\eta$ and 
$\eta^{\varepsilon}$ as the Lagrangian trajectories 
associated with the velocity solutions $u$ and $u^{\varepsilon}$.  The 
vector $\eta$ solves the first order initial value 
problem 
\begin{equation}\label{eq:eta}
    \begin{split}
        \dot{\eta} (t,x) &= u(t,\eta(t,x)), \\
        \eta(0,x) &= x,
        \end{split}
\end{equation}
where the dot represents the partial derivative with respect to time.
The vector $\eta^{\varepsilon}$ solves equation (\ref{eq:eta}) with $u$ 
replaced with the solution $u^{\varepsilon}$.   We can 
regard the differential equation (\ref{eq:eta})  as a map from the space of $u$'s (spatial or Eulerian 
description) to the space of $(\eta, \dot{\eta})$ (material or 
Lagrangian description).  The Lagrangian 
fluctuation $\xi^{\varepsilon} := \eta^{\varepsilon} \circ \eta^{-1}$ is a 
volume preserving diffeomorphism that plays the role of the Reynolds
decomposition (\ref{eq:rd}).  Since the decomposition is made in the 
Lagrangian reference frame, the interplay between the Euler and Lagrangian 
frameworks is of fundamental importance.  This interplay is 
insignificant in the Reynolds averaging 
approach where the decomposition and the averaging only occur in the Eulerian 
reference frame.  After asymptotically expanding 
$u^{\varepsilon}$ in terms of $\varepsilon$, Marsden and Shkoller take the ensemble average 
of the total kinetic energy
\[ S = \frac{1}{2} \int_0^T \int_\Omega |u^{\varepsilon}|^2 \, dx \, dt \]
over all possible solutions $u^{\varepsilon}$.  At this point, the turbulence 
closure problem is confronted.  In the Lagrangian averaging approach, the 
turbulence closure problem amounts to specifying 
the fluctuations $\xi^{\prime}$ and $\xi''$ (\, $^{\prime}$ 
denotes the derivative with respect 
to $\varepsilon$ evaluated at $\varepsilon = 0$) as functions of the mean
velocity $u$.  The generalized Taylor ``frozen turbulence'' 
hypothesis provides the solution and assumes the fluctuation 
$\xi^{\prime}$ is frozen or Lie advected in the mean flow as a 
divergence-free vector field.  In addition, the fluctuation $\xi''$ satisfies 
\[ \frac{D}{Dt} \langle \xi'' \rangle \perp u \mbox{ in } L^2(\Omega), \]
where $\langle \, \cdot \, \rangle$ is the ensemble average over the 
solutions $u^\varepsilon$ and the operator $\frac{D}{Dt}$ is the usual 
total derivative.  Marsden and Shkoller's derivation results 
in an averaged action functional which includes all terms up to order 
$\alpha^{2}$.  To conclude their calculation, the authors apply 
Hamilton's principle to the averaged action functional therefore yielding the 
LAE-$\alpha$ equations.

The derivation outlined above has recently been generalized by Bhat \emph{et. al} 
\cite{bfm04} in their computation of the 
LAE-$\alpha$ equations for a compressible fluid.  Unlike the previous 
derivation, the 
authors averaged over a tube of trajectories centered around a given 
Lagrangian flow.  The tube is constructed by specifying the Lagrangian 
fluctuation $\xi^{\varepsilon}$ at $t=0$ and deciding on a ``flow rule'' 
which evolves $\xi^{\varepsilon}$ to later time.  For example, the 
flow rule for an incompressible fluid is the frozen Taylor 
hypothesis.  In the model of compressible fluid motion, the flow rule is the 
answer to the closure problem and is chosen to be one of two different 
physical properties.  One physical property is appropriate for 
isotropic fluid conditions and the other property yields the anisotropic model 
for bounded fluid containers.  The remainder of the author's construction 
follows the outline presented in the previous paragraph.


\section{Geometric structure of the Euler and LAE-$\alpha$ 
equations.}

The Euler-Poincar\'{e} equations, developed originally by 
Poincar\'{e} \cite{p01} in his study of Euler-type equations, are determined 
once a Lagrangian map 
$L:\mathfrak{g} \rightarrow \mathbb{R}$ is specified in a Lie algebra 
$\mathfrak{g}$.  For any point $\zeta \in \mathfrak{g}$, the evolution of 
the variable $\zeta$ is determined by the Euler-Poincar\'{e} equations
\begin{equation}\label{eq:ep}
    \frac{d}{dt} \frac{\delta L}{\delta \zeta} = 
    \operatorname{ad}_{\zeta}^{*} \frac{\delta L}{\delta \zeta}.
\end{equation}
The map $\operatorname{ad}_{\zeta} : \mathfrak{g} \rightarrow 
\mathfrak{g}$, the adjoint representation of the Lie 
algebra, is the linear map $\eta \mapsto [ \zeta, \eta ]$, where 
$[ \zeta , \eta ]$ denotes the Lie bracket of $\zeta$ and $\eta$.  The map 
$\operatorname{ad}_{\zeta}^{*}$ is the dual linear map associated with 
$\operatorname{ad_{\zeta}}$.  The LAE-$\alpha$ equations and the 
classical Euler equations are the Euler-Poincar\'{e} equations posed on the 
same Lie algebra using different Lagrangian maps.  The configuration space 
that yields the appropriate Lie algebra for incompressible fluid motion was 
unknown until Arnold \cite{a66} and Ebin and Marsden \cite{em70}.

Ebin and Marsden \cite{em70} define the configuration space for 
incompressible fluid motion as
\begin{eqnarray*}
    \mathcal{D}_{\mu}^{s}(\Omega) &:=& \{
\eta(t,\cdot): \Omega \rightarrow \Omega \, | \, \eta \in H^{s}(\Omega), 
\eta \mbox{ is a bijection, }\,  \\ 
& & \eta^{-1}(t,\cdot):\Omega \rightarrow \Omega  \mbox{ is in } 
H^{s}(\Omega), \mbox{ and} \operatorname{det} D\eta =1  \}, 
\end{eqnarray*}
the group of $H^{s}$-class volume preserving diffeomorphisms on 
$\Omega$.  The group $\mathcal{D}_{\mu}^{s} (\Omega)$ (under 
composition) is a $C^{\infty}$ differentiable manifold but 
is not a Lie group \cite{em70} (right multiplication is smooth, but 
left multiplication is not).  It does however have an exponential map 
and associated Lie algebra in the usual sense of Lie groups.  The tangent 
space to $\mathcal{D}_{\mu}^{s}(\Omega)$ at the identity 
is identified with the space $\mathcal{X}_{\operatorname{div}}^{s}(\Omega)$, 
the space of $H^{s}$ divergence-free vector fields on $\Omega$ that are 
tangent to the boundary $\partial \Omega$.  
  
The Euler equations can be written as the Euler-Poincar\'{e} equations 
posed on the Lie algebra $\mathcal{X}_{\operatorname{div}}^{s} 
(\Omega)$ with the Lagrangian defined as the total kinetic energy.  
This connection was first made by Arnold \cite{a66}.  The Euler 
equations arise from an application of Hamilton's principle of least 
action to the $L^{2}$ Lagrangian
\begin{equation}\label{eq:u_lagr}
    \ell (u) = \frac{1}{2} \int_{\Omega} |u(t,x)|^{2} \, dx.
\end{equation}
Euler-Poincar\'{e} 
reduction techniques (see Marsden and Ratiu \cite{mr99}) show that 
Hamilton's principle reduces to the following variational principle with 
respect to Eulerian velocities:
\begin{equation}\label{eq:variation} 
    \delta \int_{a}^{b} \ell (u) \, dt = 0, 
\end{equation}
which should hold for all variations $\delta u$ of the form
\[ \delta u = \dot{w} + (u\cdot \nabla)w - (w \cdot \nabla)u, \]
where $w$ is a time dependent vector field representing the 
infinitesimal particle displacement vanishing at temporal endpoints.  
The Hamiltonian structure, along with references to the literature, can be 
found in Marsden and Weinstein \cite{mw83}, Arnold and Khesin 
\cite{ak98}, and Marsden and Ratiu \cite{mr99}.

Arnold \cite{a66} discovered that solutions to the Euler equation 
correspond to geodesics of $\mathcal{D}_{\mu}^{s}(\Omega)$ with respect a 
$L^{2}$ right invariant metric.  The $L^{2}$ metric is defined to be 
the weak Riemannian metric on $\mathcal{D}_{\mu}^{s}(\Omega)$ whose value at 
the identity is
\begin{equation}\label{eq:L2metric}
    \left< u(t,x), w(t,x) \right>_{L^{2}} = \int_{\Omega} u(t,x) \cdot 
    w(t,x) \, dx, \hspace{.2in} u, w \in 
    \mathcal{X}_{\operatorname{div}}^{s}(\Omega).
\end{equation}
The word weak is used here because this 
metric need not define the topology on the tangent space but may 
define a weaker topology.  The bridge between the group 
$\mathcal{D}_{\mu}^{s}(\Omega)$ and hydrodynamics is 
the following.  If $\eta(t,x) \in 
\mathcal{D}_{\mu}^{s}(\Omega)$ is a geodesic with respect to the $L^{2}$ 
right-invariant metric (\ref{eq:L2metric}) and the velocity of the fluid 
is defined as $u(t,x) = \dot{\eta}(t,\eta^{-1}(t,x))$, then velocity vector 
$u$ is a solution to the classical Euler equations
\begin{equation}\label{eq:euler}
    \begin{split}
    \frac{\partial u}{\partial t} + (u \cdot \nabla) u = - 
    \operatorname{grad} \, p + f, \\
    \operatorname{div} u = 0, \, \,\, u(0,x) = u(0),
    \end{split}
\end{equation}
with the boundary condition that $u$ is tangent to $\partial \Omega$.  
The geodesic $\eta$, which satisfies the first order differential 
equation (\ref{eq:eta}), is the particle path associated with the vector 
velocity solution to the Euler equations.  By showing the existence of a 
smooth geodesic on the group $\mathcal{D}_{\mu}^{s}(\Omega)$ for 
$s>(n/2)+1$, Ebin and Marsden \cite{em70} obtain short-time 
well-posedness of the Euler equations on a smooth 
$n$-dimensional Riemannian manifold.  

We summarize the geometric results of the Euler equations on $\Omega$ 
with the following theorem.

\begin{theorem}[see \cite{hmr98b}] \label{th:eu_geo}
    The following statements are equivalent:
\begin{quote}
    (i) The velocity vector $u\in 
    \mathcal{X}_{\operatorname{div}}^{s}(\Omega)$ solves the Euler equations 
    (\ref{eq:euler}).
    
    (ii) The particle path $\eta \in \mathcal{D}_{\mu}^{s}(\Omega)$, given 
    by (\ref{eq:eta}), is a geodesic on $\mathcal{D}_{\mu}^{s}(\Omega)$ with 
    respect to the $L^{2}$ right invariant weak metric 
    (\ref{eq:L2metric}).
    
    (iii) Hamilton's principle of least action (\ref{eq:variation}) holds 
    for variations of the form $\delta u = \dot{w} + (u \cdot \nabla) 
    w - (w \cdot \nabla)u$.
    
    (iv) The Euler-Poincar\'{e} equations
    \[ \frac{d}{dt} \frac{\delta \ell}{\delta u} = 
    \operatorname{ad}_{u}^{*} \frac{\delta \ell}{\delta u} \]
    holds for the Lagrangian $\ell$ defined by (\ref{eq:u_lagr}) on 
    the Lie algebra $\mathcal{X}_{\operatorname{div}}^{s}(\Omega)$.
\end{quote}
\end{theorem}

The Lagrangian averaged Euler equations contain a rich geometric structure 
similar to the geometric framework of the Euler equations.   The 
LAE-$\alpha$ equations were first written by Holm, Marsden, and 
Ratiu \cite{hmr98b} as an $n$-dimensional generalization of the one-dimensional 
Camassa-Holm equation.  For a given parameter $\alpha > 0$, the 
authors define the Lagrangian for the mean fluid 
velocity $u \in \mathcal{X}_{\operatorname{div}}^{s}(\Omega)$ as 
the $H^{1}$-equivalent norm
\begin{equation}\label{lae_lag}
    L(u) = \frac{1}{2} \int_{\Omega} u \cdot u + 2 \alpha^{2}
    \operatorname{Def} u : \operatorname{Def} u  \, dx, 
\end{equation}
where the differential operator $\operatorname{Def}$ is the deformation 
tensor defined by (\ref{eq:Def}) and $:$ is the contraction on two indices 
given by $a:b = a_{ij}b_{ij}$.  This $H^{1}$ Lagrangian defines a 
right-invariant weak metric on the Lie algebra of divergence-free vector spaces 
$\mathcal{X}_{\operatorname{div}}^{s}(\Omega)$.  
Holm, Marsden, and Ratiu 
\cite{hmr98b} calculate the Euler-Poincar\'{e} equations with the 
Lagrangian defined by equation (\ref{lae_lag}).  The result is the Lagrangian 
averaged Euler equations in Euclidean space
\begin{equation}\label{eq:lae}
  \begin{split}
      \frac{\partial v}{\partial t} + (u \cdot \nabla) v + & \left[ \nabla u 
      \right]^{T} \cdot v = - \operatorname{grad} \, p, \\
      v= (1&-\alpha^{2}\lap)u, \\
      \operatorname{div} \, v = 0, &\hspace{.2in} v(0,x) = v_{0},
    \end{split}
\end{equation}
with the boundary condition that $v$ is tangent to $\partial 
\Omega$.   The authors conclude that the solutions to the LAE-$\alpha$ 
equations correspond to the geodesics of 
$\mathcal{D}_{\mu}^{s}(\Omega)$ with respect to the $H^{1}$ right-invariant 
metric defined by (\ref{lae_lag}).
In \cite{s98} and \cite{s00}, Shkoller generalizes equation (\ref{eq:lae}) to a 
smooth $n$-dimensional Riemannian manifold.  By showing the existence of 
smooth geodesic flow on $\mathcal{D}_{\mu}^{s}(\Omega)$ with respect to the 
$H^{1}$ metric defined by (\ref{lae_lag}), Shkoller proves the 
smooth-in-time local (global when $n=2$)
existence and uniqueness of strong solutions with $H^{s}$ initial data, 
$s>(n/2)+1$.  


\section{Organization of the dissertation.}

In this dissertation we study the well-posedness of the Lagrangian averaged 
Navier-Stokes equation on fluid containers located in 
$\mathbb{R}^{3}$.  
There are two types of LANS-$\alpha$ equations: the anisotropic 
version in which the (fluctuation) covariance tensor is a dynamic 
variable that is coupled with the evolution equations for the mean velocity, 
and the isotropic version in which the covariance tensor is assumed 
to be a constant multiple of the identity matrix.  A brief history of the 
LANS-$\alpha$ equations is presented in Chapter \ref{chapter:iso_weak}. 

In Chapter \ref{chapter:iso_weak}, we consider the isotropic 
LANS-$\alpha$ equations on bounded domains.  We prove the global-in-time 
existence and uniqueness of weak solutions for the case of no-slip 
boundary data, extending the periodic box result of Foias, Holm, and Titi 
\cite{fht02}.  We make use of a formulation of the LANS-$\alpha$ equations on 
bounded domains given by Shkoller \cite{s00} and Marsden and 
Shkoller \cite{ms01} which reveals the correct boundary conditions.  We 
begin the proof of existence with a sequence of strong
solutions to the LANS-$\alpha$ equations.   Known to exist from 
\cite{ms01}, these vector fields trivially solve the weak formulation 
of the equations.  We use standard interpolating inequalities to 
establish estimates independent of the sequence 
index.  Classical compactness arguments enable us to conclude that the sequence 
converges in the correct space with the limit satisfying the weak 
formulation of the LANS-$\alpha$ equations.  At the end of the chapter we 
prove the existence of a nonempty, compact, convex, and 
connected global attractor.

Chapters \ref{chapter:spwp} and \ref{chapter:numerical} focus on the 
anisotropic equations.   In bounded domains, the covariance tensor 
plays a prominent role in the mechanics of fluid motion.  The 
anisotropic LANS-$\alpha$ are therefore the correct model to capture 
the large scale motion of the fluid.  There are two choices for the 
divergence-free projection of the viscosity term.  One choice is the 
classical $L^{2}$-orthogonal Leray 
projector.  In this case, Marsden and Shkoller \cite{ms03} show that 
strong solutions exist and are unique in the three-dimensional periodic box 
for a finite time interval.  In Chapter \ref{chapter:spwp}, we extend this 
result by considering the second choice of projector, the generalized Stokes 
projector.  The generalized Stokes projector, defined in detail in 
Chapter \ref{chapter:spwp}, assigns to the divergence-free vector 
field the no-slip boundary condition.  Without fluid motion on the 
boundary, the fluctuation tensor is also zero on the boundary, making the 
Stokes projector the appropriate projector for the anisotropic equations.  
The inclusion of the generalized Stokes projector viscosity term in the 
anisotropic equations on a periodic box is an essential step towards 
understanding averaged flow on bounded domains.  We prove the local-in-time 
existence and uniqueness of classical solutions to the anisotropic 
LANS-$\alpha$ on the periodic box with this 
viscosity term by using quasi-linear partial differential equation type 
methods.  We begin the proof by obtaining an approximate solution using the 
Galerkin projection of the anisotropic equations
onto a finite dimensional vector space.  The generalized Stokes projector term 
forces 
us to repose the problem in terms of the momentum rather than the velocity.  
After proving an elliptic regularity-type result, we show that 
our approximations remain in the correct vector spaces independent of 
the projection.  We use classical compactness arguments to establish 
the existence of solutions for a short period of time. 

In Chapter \ref{chapter:numerical}, we examine the anisotropic LANS-$\alpha$ 
equations in the channel and pipe geometry.  Recently, Coutand and 
Shkoller \cite{cs03a} have proposed a turbulent channel flow theory 
founded on the anisotropic equations.  By posing the 
problem in the correct functional framework, Coutand and Shkoller 
\cite{cs03a} proved that weak solutions of the anisotropic LANS-$\alpha$ 
equations exist and are unique throughout the entire channel for all time.  
In Chapter \ref{chapter:numerical}, we numerically compute strong solutions 
to the anisotropic equations by considering steady fluid flow with 
no-slip boundary conditions.  In particular, given a steady velocity vector 
that solves the Navier-Stokes equations, we numerically calculate the 
covariance tensor such that the pair solves the anisotropic LANS-$\alpha$ 
equations.  Our solutions are in good agreement with analytic results 
contained in \cite{cs03a}.  Namely, we confirm the logarithmic degeneracy rate 
of the covariance tensor near the boundary and show that in these 
elementary domains the sup-norm of the covariance tensor is unbounded in time 
near the wall.  We conclude the chapter by showing the existence of a 
shear flow velocity solution to the anisotropic LANS-$\alpha$ 
equations.  

This research was partially supported by the 
National Science Foundation under Grant~\#~DMS-0135345, KDI grant 
\# ATM-98-73133, and VIGRE Grant 0135345.

\newpage
\pagestyle{myheadings} 
\markright{  \rm \normalsize CHAPTER 2. \hspace{0.5cm} 
 Weak solutions of the LANS-$\alpha$ equations on bounded domains }
\chapter{Weak solutions of the isotropic LANS-$\alpha$ equations}
\label{chapter:iso_weak}
\thispagestyle{myheadings}
\section{Introduction.}

The isotropic Lagrangian averaged Navier-Stokes (LANS-$\a$) equations for an
incompressible viscous fluid moving in a bounded fluid container
$\Omega \subset \mathbb{R}^n$ with smooth (at least $C^3$)
boundary $\partial \Omega$ may be written as the following system
of partial differential equations:

\begin{subequations}
  \label{lans}
\begin{gather}
\partial_t u + \nabla_u u + \mathcal{U}^{\a}(u) = -
(1-\a^2 \lap)^{-1} \operatorname{grad}p - \nu Au + \mathcal{F} \,,
         \label{lans.a}\\
 \operatorname{div}u=0 \,,
         \label{lans.b}\\
u=0 \text{ on } \partial \Omega \,,
         \label{lans.c}\\
u(0,x) = u_0(x)\,,
         \label{lans.d}
\end{gather}
\end{subequations}
where
\begin{equation}\label{lans.e}
\mathcal{U}^{\a}(u) = \a^2(1-\a^2
\lap)^{-1}\operatorname{Div}(\nabla u \cdot {\nabla u}^T + \nabla
u \cdot \nabla u - {\nabla u}^T \cdot \nabla u),
\end{equation}
and $\partial_{t}$ represents the partial derivative with respect to $t$.
We use $u(t,x)$ to denote the
large-scale (or averaged) velocity field of the fluid, assumed to have
constant density.  The pressure
function $p(t,x)$  is determined (modulo constants) from the incompressibility
constraint (\ref{lans.b}).  The constant $\nu$ denotes the {\it kinematic}
viscosity of the fluid, and $\a>0$ is the spatial scale at which fluid motion
is filtered, i.e. spatial scales smaller than $\alpha$ are averaged
out.  The additional term $\mathcal{F}(x)\in
H^2 \cap H^1_0$ represents the external force acting on the system
and is assumed, for simplicity, to be time-independent.  We let
$A :=-P\lap$ denote the Stokes operator, with $P$ the
Leray projector onto divergence-free vector fields.

It has been a longstanding problem in fluid dynamics to derive a
model for the large scale motion of a fluid that averages or
course-grains the small, computationally unresolvable, scales of
the Navier-Stokes equations. The LANS-$\alpha$ equations provide
one such averaged model, and have been studied rather extensively
from both the analytical, as well numerical, points of view.  
The numerical simulations of both forced and decaying isotropic 
turbulence given by Chen {\it et al.} \cite{cfh98,
cfh99a, cfh99b, chm99} and Mohseni {\it et al.} \cite{msk03} are 
briefly reviewed in Chapter \ref{chapter:numerical}.  

The LANS-$\alpha$ equations can be written in an isotropic and 
anisotropic version.  The isotropic equations, given by equations 
(\ref{lans}) are the focus of this chapter. 
As discussed in Chapter \ref{chapter:introduction}, the inviscid ($\nu=0$) 
version of equations (\ref{lans}), known as the Lagrangian
averaged Euler (LAE-$\alpha$) or Euler-$\alpha$ equations, was first given by
Holm, Marsden, and Ratiu~\cite{hmr98b} in the case that
$\Omega= {\mathbb R}^n$ as
\begin{equation}\label{hmr}
\begin{array}{c}
\partial_t v + \nabla_u v + \left[ \nabla u \right]^T \cdot v
= -\operatorname{grad} p , \\
\operatorname{div} u = 0, \,
\end{array}
\end{equation}
where the variable
$$v = (1-\a^2 \lap) u$$
may be thought of as the momentum.
Foias, Holm, and Titi~\cite{fht02} first added viscous dissipation to
(\ref{hmr}); they argued on physical grounds that the momentum $v$ rather
than the velocity $u$, need be diffused.  By assuming periodic boundary
conditions,
they obtained the following form of the LANS-$\alpha$ equations:
\begin{equation}\label{flans}
\begin{array}{c}
\partial_t v+ \nabla_u v - \a^2 \left[\nabla u \right]^T\cdot
\lap u = -\operatorname{grad} p + \nu \lap v + g, \\
\operatorname{div} u = 0, \\
u(0,x) = u_0(x),
\end{array}
\end{equation}
with $g$ taken in $L^2$.  While yielding the correct equations on
a periodic box, the question of how to appropriately prescribe boundary data
in the no-slip $u=0$ case remained open.
Specifically, inversion of the dissipative term $\nu \lap v =
\nu (1-\alpha^2\lap) u$, a fourth-order operator, requires further
constraints than simply $u=0$ on ${\partial \Omega}$.

Shkoller \cite{s00} and Marsden and Shkoller \cite{ms01} supplied
the additional boundary condition by reformulating (\ref{flans})
as the system of equations (\ref{lans}).
In the formulation (\ref{lans}), it is clear that if $u=0$ on
${\partial \Omega}$, then
$$Au =0 \ \text{ on } \ \partial \Omega$$
as well.  This follows since each term in the inviscid equations
identically vanishes on the boundary, thanks to the inversion of
$(1-\alpha^2\lap)$ with Dirichlet boundary conditions.  The
viscous term in the formulation (\ref{lans}) was obtained by
treating the Lagrangian trajectory as a stochastic process, and
replacing deterministic time derivatives with backward-in-time
mean stochastic derivatives, exactly following the usual procedure
for obtaining the viscous dissipation term in the Navier-Stokes
equations as done by Chorin \cite{c73} and Peskin \cite{p85}.

The term ${\mathcal U}^\alpha(u)$, given in equation (\ref{lans.e}),
provides a regularization to the Navier-Stokes equations which is
dispersive, rather than dissipative, in character.
This regularization is geometric in nature, and
arises as the geodesic flow of an $H^1$ right-invariant Riemannian
metric on the Hilbert group of volume-preserving diffeomorphisms
of the fluid container (discussed in Chapter 
\ref{chapter:introduction}); as such, this
regularizer yields an a priori $L^\infty-H^1$ estimate in
three-dimensions (in general, $n$ dimensions for $n\ge 2$), and
one is thus tempted to ask whether the LANS-$\alpha$ system is
globally well-posed (even though, as is clear from (\ref{lans.a}),
no additional artificial viscosity is being added to the
Navier-Stokes equations).

In the case of the periodic box, Foias, Holm, and Titi~\cite{fht02}
proved the
global well-posedness of $H^1$ weak solutions in dimension three, but as
we noted, their formulation (\ref{flans}) did not provide the obvious extension
to bounded domains.  Using the equations (\ref{lans}), Marsden and Shkoller
\cite{ms01} proved the global well-posedness of classical solutions in
dimension three in the case of no-slip boundary  data.
In this chapter, we give an extension of that result
to the $H^1$ weak solutions of Foias, Holm, and Titi \cite{fht02}.
In particular, we prove the
global-in-time existence, uniqueness, and regularity of weak solutions
to the LANS-$\a$ equations for initial data in the class $\{ u \in
H^1_0 | \operatorname{div} u = 0 \}$.  The analogous two-dimensional
result follows trivially (as it is already known for the original
Navier-Stokes equations).

The well-posedness result leads to the existence
of a nonempty, compact, convex, and connected global $H^1$
attractor in both two- and three-dimensions.  In three-dimensions, the
global attractor has the identical bound as that obtained by
Foias, Holm, and Titi~\cite{fht02} for periodic boundary conditions.
This upper bound depends
on $1/\a$ and consequently tends to infinity as $\a$ tends to zero.
Due to the difference in the Lieb-Thirring inequality in two- and
three-dimensions, the two-dimensional bound for the global
attractor is $\a$-independent.  The global attractor of
the two-dimensional LANS-$\a$ equations is similar to the bound for the global
attractor of the two-dimensional Navier-Stokes equations.  This is shown in 
detail in Coutand \emph{et. al.} \cite{cps02}.

The remainder of this chapter is organized into two sections.  In Section 
\ref{sec:global}, we establish the global well-posedness result, while 
Section \ref{sec:attractor} is devoted to
showing the existence of a global $H^{1}$ attractor.

\section{Global well-posedness.}\label{sec:global}

In this section, we establish the global existence of unique $H^1$
weak solutions to the LANS-$\a$ equations on bounded domains
$\Omega \in \mathbb{R}^3$. Rather than using the standard Galerkin
method, we instead take a sequence of classical solutions, via the
existence result in \cite{ms01}, and prove that this sequence
converges in $C([0,T], H^1)$ to a $H^1$ weak solution of
equations (\ref{lans}) for all $T>0$.

\subsection{Notation and some classical inequalities.}

We work in the following Hilbert spaces: $\mathcal{V}^s
= H^s \cap H_0^1$ and $\mathcal{V}_{\mu}^s = \{ u \in
\mathcal{V}^s \, | \, \operatorname{div} \, u = 0 \}$ for $s \geq
1$, and $\dot{\mathcal{V}}_{\mu}^s = \{ u \in \mathcal{V}_{\mu}^s
\, | \, Au = 0 \, \mathrm{on} \, \partial \Omega \}$ for $s \geq
3$.  
We endow $\mathcal{V}^1_{\mu}$ with the following scalar
product:
\begin{equation} \label{eq:product}
    \langle f,g\rangle_{1,\a} = \int_{\Omega} f\cdot g + \a^2 \,
\mathrm{Def}\, f \,:\, \mathrm{Def} \, g \, dx, 
\end{equation}
where $\mathrm{Def} \, u$ is the rate of deformation tensor defined 
by equation (\ref{eq:Def}), $\cdot$ denotes the
usual dot product, and $:$ is the contraction of two indices, e.g.
$a : b = a_{ij} b_{ij}$.  

\ni Furthermore for any integer $s \geq 0$, we set
$$D^s u = \{D^{\beta} u : |\beta| = s \}, \hspace{.3in} ||D^s
u||_{L^p} = \sum_{|\beta| = s} ||D^{\beta} u ||_{L^p} $$ where $\beta
= (\beta_1, \dots , \beta_n)$ denotes a multi-index, and
$$|\beta| = \beta_1 + \cdots + \beta_n, \hspace{.3in} D^{\beta} =
\partial_1^{\beta_1} \cdots \partial_n^{\beta_n},$$
where 
\[ \partial_{k} := \frac{\partial}{\partial x_{k}}. \]

Throughout the dissertation, we let $C>0$ denote a generic constant.  For 
simplicity in notation, we write $u(t)=u(t,\cdot)$.  The following 
standard inequalities are used frequently and we include their definitions 
for completeness.

\noindent {\bf Gagliardo-Nirenberg inequalities~\cite{adn59} .}  Suppose
\[ \frac{1}{p} = \frac{i}{n} + a \left( \frac{1}{r} - \frac{m}{n}
\right) + (1-a)\frac{1}{q} \] where $ 1/m \leq a \leq 1$ (if
$m-n-n/r$ is an integer $\geq 1$, only $a<1$ is allowed).  Then
for $f:\Omega \rightarrow \mathbb{R}^n$,
\begin{equation}\label{gn}
||D^i f||_{L^p} \leq C ||D^m f||_{L^r}^a \cdot ||f||_{L^q}^{1-a}.
\end{equation}
\ni Two specific cases of (\ref{gn}) in
dimension three are
\begin{equation}\label{gn1}
||v||_{L^4} \leq 4 ||Dv||_{L^2}^{3/4} ||v||_{L^2}^{1/4},
\end{equation}
\begin{equation}\label{gn2}
||D^i v ||_{L^2} \leq C ||v||_{L^2}^{1-i/m} ||D^m v ||_{L^2}^{i/m}.
\end{equation}

\noindent {\bf Sobolev Embedding Theorem.} If $u\in 
H^{s}(\mathbb{R}^{n})$ for $\displaystyle{s>\frac{n}{2}+k}$, then $u\in 
C^{k}(\mathbb{R}^{n})$.  In particular for 
$\displaystyle{s>\frac{n}{2}, H^{s}(\mathbb{R}^{n})\subset 
L^{\infty}(\mathbb{R}^{n})}$ with the bound
\begin{equation}\label{eq:sobolev}
    ||u||_{L^{\infty}(\mathbb{R}^{n})} \leq C ||u||_{H^{s}(\mathbb{R}^{n})}
\end{equation}
where the constant $C$ only depends on $s$ and $n$.

\noindent {\bf Young's Inequality.} Let $1<p,q<\infty$, 
$\displaystyle{\frac{1}{p} + \frac{1}{q} =1}$.  Then 
\begin{equation}\label{eq:youngs}
    ab \leq \varepsilon a^{p} + C(\varepsilon) b^{p} \hspace{.2in} 
    (a,b>0, \varepsilon >0)
    \end{equation}
for $C(\varepsilon) = (\varepsilon p)^{-q/p}q^{-1}$.

\noindent {\bf Gronwall's Inequality.}  Let $\eta(\cdot)$ be a 
nonnegative, absolutely continuous function on $[0,T]$ which 
satisfies for a.e. $t$ the differential inequality
\[ \eta'(t) \leq \phi(t) \eta(t) + \psi(t), \]
where $\phi(t)$ and $\psi(t)$ are nonnegative, summable functions 
on $[0,T]$.  Then 
\begin{equation}\label{eq:gronwall}
 \displaystyle{   \eta (t) \leq e^{\int_{0}^{t} \phi(s)\, ds} \left[ \eta(0) + 
    \int_0^{t} \psi (s) \, ds \right]}
\end{equation}
for all $0 \leq t \leq T$.

\subsection{Three equivalent forms of the LANS-$\a$
equations.}\label{lansform} 

Three equivalent forms of the LANS-$\a$
equations will be useful to us.

\ni {\bf LANS-1:}
\begin{equation}\label{lans1}
\begin{array}{c}
\partial_t u + \nabla_u u + \mathcal{U}^{\a}(u) = -\nu  Au -
(1-\a^2 \lap)^{-1} \operatorname{grad} p +(1-\a^2 \lap)^{-1}f\\
\operatorname{div}u(t,x)=0, \\
u = 0 \text{ on } \partial \Omega, \\
u(0,x) = u_0(x),
\end{array}
\end{equation}
where $f\in L^2$.  (It is convenient to replace $\mathcal{F}$ in equation
(\ref{lans.a}) with $(1-\alpha^2\lap)^{-1}f$; there is no loss in
generality as $(1-\alpha^2\lap)$ with domain $H^2\cap H^1_0$ is
an isomorphism.)

The Stokes operator $Au = -P\lap u$, is the Leray projection of
$-\lap u$ onto divergence-free vector fields and has domain $D(A)
= H^2 \cap H^1_0$.  As we noted above, when $u=0$ on ${\partial \Omega}$,
then $Au$ must also equal zero on $\partial \Omega$.

\ni {\bf LANS-2:} This form is equivalent to LANS-1 in view of our
remark that LANS-1 implies $Au =0$ on $\partial \Omega$:

\begin{equation}\label{lans2}
\partial_t(1-\a^2 \lap)u + \nabla_u \left[(1-\a^2 \lap)u \right] - \a^2
\nabla u^T \cdot \lap u
=\nu (1-\a^2 \lap)Au - \operatorname{grad} p + f ,
\end{equation}
together with the constraint $\operatorname{div}u(t,x) = 0$ and
boundary data $ u = Au = 0$ on $\partial \Omega$.
Note that when the domain $\Omega$ is the period box ${\mathbb T}^3$,
the Stokes operator is given by $-\lap$, and formulation (\ref{lans2}) 
reduces to equations (\ref{flans}), the LANS-$\alpha$ equations used by Foias, 
Holm, and Titi \cite{fht02}.

\ni {\bf LANS-3:} This form is the analog of the Navier-Stokes
equations written in terms of the Helmholtz-Hodge projection:
 \begin{equation}\label{lans3}
\partial_t u + \nu Au + \mathcal{P}^{\a}\left[\nabla_u u + \mathcal{U}^{\a}(u) -
(1-\a^2 \lap)^{-1}f \right] = 0,
\end{equation}
where for $s\geq 1$, $\mathcal{P}^{\a}$ is the Stokes projector
defined below.
\begin{definition}\label{def:stokes} For $s\ge 1$,  we
let $\mathcal{P}^\alpha:\mathcal{V}^s \rightarrow
\mathcal{V}^s_{\mu}$ denote the {\it Stokes projector}, a continuous
$\langle \cdot, \cdot \rangle_{1,\alpha}$-orthogonal idempotent
operator (see Proposition 1 of \cite{s00}).  It is defined as
\[ \mathcal{P}^{\a}(w) = w - (1-\a^2 \lap)^{-1} \operatorname{grad} p =v,\]
where $(v,p)$ solve the {\it Stokes problem}: given $w\in
\mathcal{V}^s$, there is a unique vector field $v\in
\mathcal{V}^s_{\mu}$ and a function $p$ (unique up to an additive
constant) such that
\begin{equation}
\nonumber
\begin{array}{c}
(1-\a^2\lap)v + \operatorname{grad} p = (1-\a^2 \lap)w, \\
\operatorname{div} v=0,\\
v=0 \ \operatorname{on}  \ \partial\Omega.
\end{array}
\end{equation}
\end{definition}

\subsection{Results.}
We begin with two elementary lemmas. 
\begin{lemma}\label{u0conv}
$\dot{\mathcal{V}}_{\mu}^4$ is dense in $\mathcal{V}^1_{\mu}$.
\end{lemma}
\noindent {\bf Proof:} Let $v \in
\mathcal{V}^1_{\mu}$.  We find a sequence $v_m \in
\dot{\mathcal{V}}^4_{\mu}$ which converges to $v$ in $H^1$. The
proof makes use of the Lax-Milgram Theorem to provide a compact
operator from $\mathcal{V}^1_{\mu}$ to
$\dot{\mathcal{V}}^4_{\mu}$.  Define the bilinear form
$E:\dot{\mathcal{V}}_{\mu}^4 \times \dot{\mathcal{V}}_{\mu}^4
\rightarrow \mathbb{R}$ by
$$E[u,v] = \langle u , v \rangle_{H^4}.$$
The function space $\dot{\mathcal{V}}_{\mu}^4$ is a closed
subspace of $H^4$ and therefore a Hilbert space endowed with the
usual $H^4$ topology.  Since $E$ is defined as the usual inner
product on $H^4$, it is coercive and continuous.  We define, using
the Riesz-Representation Theorem, the bounded linear function
$\tilde{f} \in \left( \dot{\mathcal{V}}^4_{\mu} \right)^{\ast}$ by
$$\tilde{f}(v) = ( f,v )_{L^2}  \hspace{.2in} \forall v
\in \dot{\mathcal{V}}_{\mu}^4$$ where $f\in
\dot{\mathcal{V}}_{\mu}^4$.  By the Lax-Milgram Theorem,  there
exists a unique $u(f) \in \dot{\mathcal{V}}_{\mu}^4$ such that
$$E[u,v] = \langle u, v \rangle_{H^4} = ( f, v )_{L^2},
\hspace{.2in} \forall v\in \dot{\mathcal{V}}_{\mu}^4.$$ This is
the weak formulation of an elliptic problem $Lu=f$ where $L$ is an
eighth order elliptic differential operator.  The map $L^{-1}:f
\mapsto u$ is a continuous map from $\dot{\mathcal{V}}_{\mu}^4
\rightarrow \dot{\mathcal{V}}_{\mu}^4$.  Since $H^4$ is compactly
embedded in $H^1$, $L^{-1}$ is a compact operator from
$\mathcal{V}^1_{\mu} \rightarrow \dot{\mathcal{V}}^4_{\mu} $, and
consequently has a discrete spectrum of eigenfunctions $\{e_i\}_{i
\in \mathbb{N}}$.  The eigenfunctions $\{e_i\}_{i \in \mathbb{N}}$
form a Hilbert basis of $\mathcal{V}^1_{\mu}$ and therefore we
write $v=\displaystyle{ \sum_{i=1}^{\infty} c_i e_i}$.  The
ellipticity of $L^{-1}$ implies that $\{e_i\}_{i \in \mathbb{N}}$
are also eigenfunctions of $L$, allowing us to conclude that
$\{e_i\}_{i \in \mathbb{N}} \in \dot{\mathcal{V}}^4_{\mu}$.
Consequently, the sequence $v_m = \displaystyle{ \sum_{i=1}^{m}
c_i e_i} \in \dot{\mathcal{V}}^4_{\mu}$ converges in $H^1$ to $v$
as $m \rightarrow \infty$. \qed

\begin{lemma}\label{lem1} Let $\lambda_1 = \displaystyle{ \inf_{v \in
\mathcal{V}^1_{\mu}}
\frac{ ||\nabla v||_{L^2}^2 }{||v||_{L^2}^2 } }$ be the smallest eigenvalue of
the Stokes operator.  Then for all $v \in \mathcal{V}^2_{\mu}$,
\begin{equation*}
||\nabla v||_{L^2}^2 + \a^2 ||Av||_{L^2}^2 \geq \lambda_1 \left\{ 
||v||_{L^2}^2 +\a^2 ||\nabla v||_{L^2}^2 \right\}
\end{equation*}
\end{lemma}
 
\noindent {\bf Proof:} We just need to show
that $||Av||^2_{L^2} \geq \lambda_1 ||\nabla v||_{L^2}^2$ (since
$||\nabla v||_{L^2}^2 \geq \lambda_1 ||v||_{L^2}^2$ by definition).
Let $\{ h_i \}_{i\in \mathbb{N}}$ be eigenfunctions associated to
the Stokes operator $A$ with corresponding eigenvalues
$\lambda_i$.  For all $v \in \mathcal{V}^2_{\mu}$, we write
$\displaystyle{ v= \sum_{i=1}^{\infty} c_i h_i }$.  Then $Av
=\displaystyle{ \sum_{i=1}^{\infty} c_i \lambda_i h_i}$. The
orthogonality property of $h_i$ implies $(Av,v) = \displaystyle{
\sum_{i=1}^{\infty} c_i^2 \lambda_i = ||\nabla v||_{L^2}^2 }$ and
consequently,
\[ ||Av||_{L^2}^2 =  \sum_{i=1}^{\infty} c_i^2 \lambda_i^2 \geq  \left(\sum_
{i=1}^{\infty} c_i^2 \lambda_i \right) \lambda_1 \geq \lambda_1
(Av,v) = \lambda_1 ||\nabla v||_{L^2}^2 \] which completes the
proof. \qed

\noindent Weak solutions to the isotropic LANS-$\alpha$ equations 
(\ref{lans}) are defined below.

\begin{definition}
Let $f\in L^2(\Omega)$ and $u_0 \in \mathcal{V}^1_{\mu}$. For any
$T>0$, a function $$u\in C([0,T];\mathcal{V}^1_{\mu}) \cap
L^2([0,T];D(A))$$ with $\frac{du}{dt} \in L^2([0,T]; L^2) \cap
L^{\infty}((0,T);\mathcal{V}^1_{\mu})$ and $u(0)=u_0$ is said
to be a {\it weak solution to the LANS-$\a$ equations} with
initial data $u_0$ in the interval $[0,T]$ provided

\begin{eqnarray}\label{weak1}
\left\langle \frac{d}{dt}(1-\a^2 \lap)u, w \right\rangle_{D(A)}  &+& \nu
\left\langle (1-\a^2 \lap)Au, w\right\rangle_{D(A)} \nonumber \\
&+&  \left\langle B(u,(1-\a^2\lap)u), w\right\rangle_{D(A)}  =  (f,w)
\end{eqnarray}

\ni for every $w=w(x) \in D(A)$, the domain of the Stokes operator, and for 
almost every $t \in [0,T]$ with
$$B(u,v) = \nabla_u v + \nabla u^T \cdot v;$$ moreover, $u(0) =
u_0$ in $\mathcal{V}^1_{\mu}$. Here, the equation (\ref{weak1}) is
understood in the following sense: for every $t_0, t \in [0,T]$,
\begin{eqnarray*}
\displaystyle{ \left\langle u(t), (1-\a^2 \lap)w \right\rangle}
_{D(A)} & + & \displaystyle{ \int_{t_0}^t \left\langle
B(u(s),(1-\a^2\lap)
u(s)), w\right\rangle_{D(A)} \, ds } \\
 & = & \displaystyle{ -\int_{t_0}^t \nu \left\langle Au(s), (1-\a^2 \lap)w
\right\rangle_{D(A)} \, ds + \int_{t_0}^t (f,w) \, ds. }
\end{eqnarray*}
\end{definition}

\noindent For our proof, we shall make use of the following results from Marsden 
and Shkoller \cite{ms01}
\begin{lemma}[Theorem 5.2 of~\cite{ms01}] \label{classical}
For $u_0 \in \dot{\mathcal{V}}^s_{\mu}$, $s\in [3,5)$, and $f \in
L^2(\Omega)$, there exists a unique solution $u$ to equation
(\ref{lans1}) in $C([0,\infty), \dot{\mathcal{V}}^s_{\mu})$.
\end{lemma}
and
\begin{lemma}[Lemma 5.1 of~\cite{ms01}] \label{Ulemma}
For $s \geq 3$, $\mathcal{U}^{\alpha} : \mathcal{V}^{s} \rightarrow 
\mathcal{V}^{s}$ and $\mathcal{U}^{\alpha} : \mathcal{V}^{2} 
\rightarrow H^{1+\sigma}$ for $\sigma \in \left( 0, \frac{1}{3} 
\right)$. 
\end{lemma}

\noindent We can now state our main result.

\begin{theorem}\label{weakth}
For $f \in L^2(\Omega)$ and $u_0 \in \mathcal{V}^1_{\mu}$, there
exists a unique weak solution
$$u \in C([0,\infty),\mathcal{V}^1_{\mu}) \cap L^{\infty} ((s,\infty),
\dot{\mathcal{V}}_{\mu}^3), \hspace{.3in} \forall s>0$$ to
equation (\ref{lans1}). The solution depends continuously on the
initial data $u_0$.
\end{theorem}

\ni {\bf Proof:} Consider a sequence of initial velocity 
$\{u^{\epsilon}_0\}_{\epsilon=1}^{\infty} \in
\dot{\mathcal{V}}_{\mu}^4$ and force
$\{f^{\epsilon}\}_{\epsilon=1}^{\infty} \in H^1$ such that
$u_0^{\epsilon} \rightarrow u_0$ (by Lemma \ref{u0conv}) and
$f^{\epsilon}
 \rightarrow f$.  By Lemma \ref{classical},
for each $\epsilon \in \mathbb{N}$ there exists a unique solution
$u^{\epsilon}$ of the LANS-$\a$ equations in $C([0,\infty),
\dot{\mathcal{V}}_{\mu}^4)$.  This solution has sufficient regularity to 
satisfy the weak formulation:
\begin{eqnarray*}
\displaystyle{ \left\langle u^{\epsilon}(t), (1-\a^2 \lap)w
\right\rangle_{D(A)}} & + & \displaystyle{ \int_{t_0}^t
\left\langle B(u^{\epsilon}(s), (1-\a^2\lap)u^{\epsilon}(s)),
w\right\rangle_{D(A)} \, ds } \\
 & = & \displaystyle{ -\int_{t_0}^t \nu \left\langle Au^{\epsilon}(s),
(1-\a^2 \lap)
w\right\rangle_{D(A)} \, ds + \int_{t_0}^t (f^{\epsilon},w) \, ds, }
\end{eqnarray*}
for all $w \in D(A)$.  In order to show the limit as $\epsilon \rightarrow
\infty$ is also a weak
solution, we need to develop the appropriate energy estimates.
\smallskip

\ni {\bf An $H^1$ Estimate.}  Define $a^{\epsilon} =
\nabla_{u^{\epsilon}} u^{\epsilon} +
\mathcal{U}^{\a}(u^{\epsilon}) - (1-\a^2 \lap )^{-1}
f^{\epsilon}$.  Since $u^{\epsilon} \in H^{1}_{0}$ we may use it as a 
test function. Taking the $\langle \cdot, \cdot\rangle_{1,\a}$
inner product, defined by (\ref{eq:product}), of LANS-3 together with 
$u^{\epsilon}$,
\begin{equation*}
\left\langle \frac{d}{dt}u^{\epsilon}, u^{\epsilon}
\right\rangle_{1,\a} + \nu \langle Au^{\epsilon},
u^{\epsilon} \rangle_{1, \a } + \langle
\mathcal{P}^{\a}(a^{\epsilon}), u^{\epsilon} \rangle_{1,\a} = 0.
\end{equation*}
By the definition of $\langle \cdot, \cdot \rangle_{1,\a}$,
\begin{equation*}
\left\langle \frac{d}{dt}u^{\epsilon}, u^{\epsilon}
\right\rangle_{1,\a} + \nu \langle Au^{\epsilon},
u^{\epsilon}\rangle_{1,\a} = \frac{1}{2} \frac{d}{dt} \left[
||u^{\epsilon}||_{L^2}^2 + \a^2 ||\nabla u^{\epsilon}||_{L^2}^2
\right] + \nu \left[ ||\nabla u^{\epsilon}||_{L^2}^2 +\a^2
||Au^{\epsilon}||_{L^2}^2 \right].
\end{equation*}
Notice that the previous integration by parts used the additional
boundary condition $Au^{\epsilon}=0$. Using the properties of
$\mathcal{P}^{\a}$, we find
$$\langle \mathcal{P}^{\a}( a^{\epsilon}) , u^{\epsilon}
\rangle_{1, \a} = - (f^{\epsilon}, u^{\epsilon}).$$ Therefore,
\begin{equation}\label{eq1}
\displaystyle{ \frac{1}{2} \frac{d}{dt} \Big[ || u^{\epsilon} ||_{L^2}^2 +
\a^2 || \nabla u^{\epsilon} ||_{L^2}^2 \Big] + \nu \Big[ ||\nabla
u^{\epsilon}||_{L^2}^2 + \a^2 || Au^{\epsilon}||_{L^2}^2
\Big] = (f^{\epsilon}, u^{\epsilon}) }.
\end{equation}
By Poincar\'{e}'s inequality and Young's inequality (\ref{eq:youngs}),
\begin{equation}\label{eq2}
\displaystyle{ \frac{d}{dt} \Big[ || u^{\epsilon} ||_{L^2}^2 + \a^2 ||
\nabla u^{\epsilon} ||_{L^2}^2 \Big] + \nu \Big[ ||\nabla 
u^{\epsilon} ||_{L^2}^2
 + \a^2 ||Au^{\epsilon}||_{L^2}^2 \Big] \leq K_1 }
\end{equation}
where $K_1 = \frac{||f^{\epsilon}||_{L^2}^2}{ \nu {\lambda_1}}$.  By Lemma
\ref{lem1},
\[ \displaystyle{ \frac{d}{dt} \Big[ || u^{\epsilon} ||_{L^2}^2 + \a^2 ||
\nabla u^{\epsilon}
||_{L^2}^2 \Big] + \nu \lambda_1 \Big[ ||u^{\epsilon}||_{L^2}^2 + \a^2 || 
\nabla u^{\epsilon}||_{L^2}^2 \Big] \leq K_1 }. \]

\ni Using Gronwall's inequality (\ref{eq:gronwall}), we get the $H^1$ estimate
\begin{eqnarray}\label{eq3}
\displaystyle{ || u^{\epsilon}(t) ||_{L^2}^2 + \a^2 || \nabla u^{\epsilon}(t)
||_{L^2}^2 } &\leq& \displaystyle{ e^{-\nu {\lambda}_1 t}\left(
||u^{\epsilon}(0)||_{L^2}^2 + \a_2||\nabla u^{\epsilon}(0) ||_{L^2}^2 \right)
+ \frac{ K_1}{\nu {\lambda}_1}(1-e^{-\nu {\lambda}_1
t})} \nonumber \\
& \leq & \displaystyle{ k_1 :=
||u^{\epsilon}_0||_{L^2}^2+\a^2||\nabla u^{\epsilon}(0)||_{L^2}^2 +
\frac{ K_1}{\nu {\lambda}_1}}
\end{eqnarray} where ${\lambda}_1$ is defined in Lemma \ref{lem1}.
Therefore we have proved that
$$u^{\epsilon} \in L^{\infty}([0,\infty), \mathcal{V}^1_{\mu})$$
independently of $\epsilon$.

\ni {\bf An $H^2$ Estimate.}  Since $Au^{\epsilon}$ is
divergence-free, we take the $L^2$ inner product of LANS-2
together with $Au^{\epsilon}$ to get
\begin{eqnarray}\label{h2}
\frac{1}{2} \frac{d}{dt} \Big[ ||\nabla u^{\epsilon}||_{L^2}^2  &+&  \a^2
||Au^{\epsilon}||_{L^2}^2 \Big] + \nu \left[ ||Au^{\epsilon}||_{L^2}^2 + \a^2
||\nabla Au^{\epsilon}||_{L^2}^2 \right] \nonumber\\ & &+ (B(u^{\epsilon},
(1-\a^2 \lap)u^{\epsilon}), Au^{\epsilon}) = (f^{\epsilon}, Au^{\epsilon}).
\end{eqnarray}
By estimate (\ref{eq3}), we use the Sobolev embedding result 
$W^{1,\infty} \subset L^{4}$ to bound the nonlinear term of equation 
(\ref{h2})
\begin{eqnarray}\label{eq4}
|(\nabla_{u^{\epsilon}}u^{\epsilon}, Au^{\epsilon})| &\leq& ||\nabla
u^{\epsilon}||_{L^2}||u^{\epsilon}||_{L^4}||Au^{\epsilon}||_{L^4}
\leq Ck_1||u^{\epsilon}||_{H^1}||Au^{\epsilon}||_{L^4} \nonumber \\
&\leq& Ck_1^2 ||\nabla Au^{\epsilon}||^{3/4}_{L^2}
||Au^{\epsilon}||^{1/4}_{L^2} \nonumber \\
&\leq& \frac{\nu \a^2}{4}||\nabla Au^{\epsilon}||_{L^2}^2 +
\frac{\nu}{4}||Au^{\epsilon}||_{L^2}^2 + C(||u^{\epsilon}||_{H^1}),
\end{eqnarray}
may be estimated by the Gagliardo-Nirenberg inequality (\ref{gn}) and a repeated 
use of Young's inequality.  The boundary conditions, the
incompressibility constraint, and one integration by
parts lead to
\begin{equation}\label{eq5}
|(\nabla_{u^{\epsilon}}\lap u^{\epsilon}, Au^{\epsilon})| = 0.
\end{equation}
For the third term, using (\ref{gn}), (\ref{gn1}), and (\ref{eq3}) we have the
estimate,
\begin{eqnarray}\label{eq6}
|\left( ({\nabla u^{\epsilon}}^T) \cdot \lap u^{\epsilon}, Au^{\epsilon} 
\right)|
&\leq& ||\nabla \lap u^{\epsilon}||_{L^2}||u^{\epsilon}||_{L^4}
||Au^{\epsilon}||_{L^4} \leq C||u^{\epsilon}||_{H^1} ||\nabla
Au^{\epsilon}||^{7/4}_{L^2} ||Au^{\epsilon}||^{1/4}_{L^2} \nonumber \\
&\leq& Ck_1  ||\nabla Au^{\epsilon}||^{23/12}_{L^2}
||u^{\epsilon}||^{1/12}_{L^2} \leq C k_1^{13/12} ||\nabla
Au^{\epsilon}||^{23/12}_{L^2} \nonumber \\
&\leq& \frac{\nu}{4}||\nabla Au^{\epsilon}||^2_{L^2} +
C(||u^{\epsilon}||_{H^1}).
\end{eqnarray}
Again Young's inequality implies
\[ |(f^{\epsilon}, Au^{\epsilon})| \leq ||f^{\epsilon}||_{L^2}||
Au^{\epsilon}||_{L^2} \leq \frac{1}{\nu}||f^{\epsilon}||_{L^2}^2 +
\frac{\nu}{4} ||Au^{\epsilon}||_{L^2}^2. \] Therefore by
(\ref{eq4})-(\ref{eq6}), the equation (\ref{h2}) becomes
\begin{eqnarray*}
\frac{1}{2} \frac{d}{dt} \Big[ ||\nabla u^{\epsilon}||_{L^2}^2 &+& \a^2 ||
Au^{\epsilon}||_{L^2}^2
\Big] + \nu \left[ ||Au^{\epsilon}||_{L^2}^2 + \a^2 ||\nabla
Au^{\epsilon}||_{L^2}^2 \right] \\
&\leq& |(B(u^{\epsilon}, (1-\a^2 \lap)u^{\epsilon}), Au^{\epsilon})| +
|(f^{\epsilon}, Au^{\epsilon})| \\
&\leq&  \frac{\nu \a^2}{2}||\nabla Au^{\epsilon}||_{L^2}^2 + \frac{\nu}{2}
||Au^{\epsilon}||_{L^2}^2 + \frac{1}{\nu}||f^{\epsilon}||_{L^2}^2 + C.
\end{eqnarray*}
Hence by Lemma \ref{lem1} we have
\begin{eqnarray}\label{eq7}
\frac{d}{dt} \Big[ ||\nabla u^{\epsilon}||_{L^2}^2 + \a^2||A
u^{\epsilon}||_{L^2}^2 \Big] + \nu \lambda_1 \Big[||\nabla
u^{\epsilon}||_{L^2}^2 + \a^2 ||Au^{\epsilon}||_{L^2}^2\Big] \leq
K_2, \\
\frac{d}{dt} \left[ \left( ||\nabla u(t)||^2_{L^2} + \a^2
||Au(t)||^2_{L^2} \right)e^{\nu \lambda_1 t} \right] \leq K_2 e^{\nu
\lambda_1 t}. \nonumber
\end{eqnarray}
where $K_2 = \frac{2}{\nu}||f^{\epsilon}||^2_{L^2}+C$.
Since$||Au||^2_{L^2}$ is not necessarily bounded at $t=0$, for every
$s_1>0$ we integrate both sides from $[s_1, t]$,
\begin{equation}\label{h2bd}
||\nabla u(t)||^2_{L^2} + \a^2 ||Au(t)||^2_{L^2} \leq \left( ||\nabla
u(s_1)||^2_{L^2} + \a^2 ||Au(s_1)||^2_{L^2}\right) e^{\nu \lambda_1
(s_1 - t)} + \nu \lambda_1 K_2.
\end{equation}

\ni Thus, independently of $\epsilon$,
$$u^{\epsilon} \in L^{\infty}((s_1, \infty),
\mathcal{V}^2_{\mu}(\Omega)),$$ for all $s_1>0$

\ni {\bf An $H^3$ Estimate.} Since $A^2u^{\epsilon}$ is not
necessarily equal to zero on $\partial \Omega$, we do not use it
to derive an $H^3$-estimate. To achieve this estimate we make use
of the Ladyzhenskaya method \cite{l63}.  Let $u_t$ denote
$\partial_t u$, and differentiate the equations LANS-2 with
respect to time to get
\begin{eqnarray*}
\partial_t(1-\a^2 \lap)u^{\epsilon}_t &+& \nabla_{u^{\epsilon}_t}(1-\a^2
\lap)u^{\epsilon} +
\nabla_{u^{\epsilon}}(1-\a^2 \lap)u_t^{\epsilon} - \a^2 (\nabla
{u^{\epsilon}_t})^T \cdot \lap u^{\epsilon} \\
&-& \a^2 (\nabla {u^{\epsilon}_t})^T \cdot \lap u^{\epsilon}_t = -
\mathrm{grad} \, p_t - \nu (1-\a^2 \lap)Au_t^{\epsilon}.
\end{eqnarray*}
Noting that $u_t^{\epsilon} \in D(A)$, we take the $L^2$ inner product with
$u_t$ to get
\begin{eqnarray}\label{h3eq0}
\frac{1}{2}\frac{d}{dt}\left[ ||u^{\epsilon}_t||_{L^2}^2 \right. &+& \left. \a^2
||\nabla u^{\epsilon}_t||_{L^2}^2 \right] + \nu \left[||\nabla
u^{\epsilon}_t||_{L^2}^2 + \a^2 || Au_t^{\epsilon}||_{L^2}^2 \right] \nonumber \\
&\leq& |(u^{\epsilon}_t, \nabla_{u_t^{\epsilon}} u^{\epsilon} )| + \a^2 |(\lap
u_t^{\epsilon}, \nabla_{u^{\epsilon}} u^{\epsilon}_t +
\nabla_{u^{\epsilon}_t} u^{\epsilon})|.
\end{eqnarray}
Before estimating the right hand side, we would like a bound for
$||u_t^{\epsilon}||_{L^2}^2$ which will be used in the later
computations.  For each $\epsilon$, $u^{\epsilon}$ satisfies the
equation LANS-3 and consequently,
\begin{eqnarray}\label{h3eq1}
||u_t^{\epsilon}||_{L^2}^2 &\leq& \nu
|(Au^{\epsilon},u^{\epsilon}_t)| + |(\nabla_{u^{\epsilon}}
u^{\epsilon}, u^{\epsilon}_t)|
\nonumber \\
& & + |(\mathcal{U}^{\a}(u^{\epsilon}), u^{\epsilon}_t)| +
|((1-\a^2\lap)^{-1}f^{\epsilon}, u^{\epsilon}_t)|.
\end{eqnarray}
The first term is simply $\nu |(Au^{\epsilon},u^{\epsilon}_t)|
\leq \frac{1}{8}||u^{\epsilon}_t||_{L^2}^2 + 2\nu^2
||Au^{\epsilon}||_{L^2}^2 $ by Young's inequality.  Using
(\ref{gn1}) and Young's inequality, the second term becomes,
\begin{eqnarray*}
|(\nabla_{u^{\epsilon}} u^{\epsilon}, u^{\epsilon}_t)| &\leq& 16
||\nabla u^{\epsilon}||_{L^2} ||u^{\epsilon}||^{1/4}_{L^2}
||u^{\epsilon}_t||_{L^2}||Au^{\epsilon}||^{3/4}_{L^2} \\ &\leq& 
\nu^2 ||Au^{\epsilon}||_{L^2}^2
 + C ||\nabla u^{\epsilon}||^{8/5}_{L^2} ||u^{\epsilon}||^{2/5}_{L^2} 
 ||u^{\epsilon}_t||^{8/5}_{L^2} \\
&\leq& \nu^2||Au^{\epsilon}||_{L^2}^2
+\frac{1}{8}||u^{\epsilon}_t||_{L^2}^2 + Ck_1^{10},
\end{eqnarray*}
where $k_1$ is the time independent $H^1$ bound (\ref{eq3}). Let
$v^{\epsilon}=(1-\a^2 \lap)^{-1}u^{\epsilon}_t$. Then, by
integration by parts, the third term of (\ref{h3eq1}) is
\begin{eqnarray*}
|(\mathcal{U}^{\a}(u^{\epsilon}), u^{\epsilon}_t)| &\leq& \a^2
|(\lap u^{\epsilon} \cdot \nabla u^{\epsilon}, v^{\epsilon})|
 \leq \a^2 ||u^{\epsilon}||_{H^2}^2  ||v^{\epsilon}||_{H^1} \\
&\leq& \frac{\a^2}{4} ||v^{\epsilon}||_{H^1}^2 + \a^2||u^{\epsilon}||_{H^2}^4.
\end{eqnarray*}
By definition, $v^{\epsilon}-\a^2 \lap v^{\epsilon} =
u^{\epsilon}_t$. Taking the $L^2$ inner product of both sides with
$v^{\epsilon}$ we have an $H^1$ estimate
\[ ||v^{\epsilon}||_{L^2}^2 + \a^2 ||\nabla v^{\epsilon}||_{L^2}^2 =
(u^{\epsilon}_t, v^{\epsilon}) \leq \frac{1}{2}
||u^{\epsilon}_t||_{L^2}^2 + \frac{1}{2} ||v^{\epsilon}||_{L^2}^2. \]
Hence $||v^{\epsilon}||_{H^1}^2 \leq \frac{1}{2\a^2} 
||u^{\epsilon}_t||_{L^2}^2$ and therefore
\[ |(\mathcal{U}^{\a}(u^{\epsilon}), u^{\epsilon}_t)| \leq \frac{1}{8}
||u^{\epsilon}_t||_{L^2}^2 + \a^2||u^{\epsilon}||_{H^2}^4. \]
The last term of (\ref{h3eq1}) can be estimated in a similar way
\begin{equation*}
|((1-\a^2\lap)^{-1}f^{\epsilon}, u^{\epsilon}_t)| =
|(f^{\epsilon}, v^{\epsilon})| \leq   \frac{1}{8}||v^{\epsilon}||
_{L^2}^2 + 2||f^{\epsilon}||_{L^2}^2 \leq
\frac{1}{8}||u^{\epsilon}_t||_{L^2}^2 + 2||f^{\epsilon}||_{L^2}^2.
\end{equation*}
Combining the above estimates, we obtain the bound
\begin{equation}\label{utbnd}
||u^{\epsilon}_t||_{L^2}^2 \leq 6\nu^2 ||Au^{\epsilon}||_{L^2}^2 +
2\a^2||u^{\epsilon}||_{H^2}^2 + 4||f^{\epsilon}||_{L^2}^2 + Ck_1^{10},
\end{equation}
From the $H^2$ estimate (\ref{h2bd}) we conclude that
$u_t^{\epsilon} \in L^{\infty} \left( (s_1, \infty), L^2 \right)$
for all $s_1>0$. Using this result, we now estimate each of
the terms on the right hand side of (\ref{h3eq0}).  Using
(\ref{gn})-(\ref{gn2}) and Young's inequality, we get
\begin{eqnarray*}
|(u^{\epsilon}_t, \nabla_{u_t^{\epsilon}} u^{\epsilon} )| &+& \a^2
|(\lap u_t^{\epsilon}, \nabla_{u^{\epsilon}} u^{\epsilon}_t +
\nabla_{u^{\epsilon}_t} u^{\epsilon})|  \\
&\leq& C \left [ ||\nabla u_t^{\epsilon}||^{3/4}_{L^2}
||u_t^{\epsilon}||^{5/4}_{L^2} ||u^{\epsilon}||_{H^2} +  \a^2
||Au_t^{\epsilon}||^{15/8}_{L^2} ||u_t^{\epsilon}||^{1/8}_{L^2}
||u^{\epsilon}||_{H^1} \right. \\
&+& \left. \a^2 ||Au_t^{\epsilon}||_{L^2} ||\nabla
u_t^{\epsilon}||^{3/4}_{L^2} ||u^{\epsilon}||^{1/4}_{L^2}
||u^{\epsilon}||_{H^1} \right] \\
&\leq& \frac{\nu}{2} ||\nabla u_t^{\epsilon}||^2_{L^2} + \frac{\nu
\a^2}{2}||Au_t^{\epsilon}||^2_{L^2} + C_{s_1},
\end{eqnarray*}
where $C_{s_1}$ is a constant multiple of the $H^2$ bound
depending only on the lower time bound $s_1>0$.  Therefore
(\ref{h3eq0}) becomes
\begin{equation*}
\frac{d}{dt}\left[ ||u^{\epsilon}_t||_{L^2}^2 + \a^2
 ||\nabla u^{\epsilon}_t||_{L^2}^2 \right]
+ \nu \left[||\nabla u^{\epsilon}_t||_{L^2}^2 + \a^2 ||
Au_t^{\epsilon}||_{L^2}^2
 \right] \leq C_{s_1}.
\end{equation*}
By Lemma \ref{lem1},
\begin{eqnarray*}
\frac{d}{dt}\left[ ||u^{\epsilon}_t||_{L^2}^2 + \a^2 ||\nabla
u^{\epsilon}_t||_{L^2}^2 \right] + \nu \lambda_1 \left[
||u^{\epsilon}_t||_{L^2}^2 + \a^2 ||
 \nabla u_t^{\epsilon}||_{L^2}^2 \right]
\leq C_{s_1} \\
\frac{d}{dt} \left[ \left(||u_t^{\epsilon}||^2_{L^2} + \a^2 ||\nabla
u_t^{\epsilon}||_{L^2}^2 \right) e^{\nu \lambda_1 t} \right] \leq
C_{s_1} e^{\nu \lambda_1 t}.
\end{eqnarray*}
For every $s_2>0$, we integrate from $[s_2, t]$ and obtain the
$H^1$ bound of $u_t^{\epsilon}$,
\begin{equation*}
||u_t^{\epsilon}(t)||^2_{L^2} + \a^2 || \nabla
u_t^{\epsilon}(t)||^2_{L^2} \leq \left(
||u_t^{\epsilon}(s_2)||^2_{L^2} + \a^2 ||\nabla
u_t^{\epsilon}(s_2)||^2_{L^2} \right) e^{\nu \lambda_1 (s_2-t)} +
\nu \lambda_1 C_{s_1}
\end{equation*}
Hence, independently of $\epsilon$, $$u^{\epsilon}_t \in
L^{\infty}((s_2, \infty), \mathcal{V}_{\mu}^1),$$ for all $s_2>0$.

Let $s = \operatorname{min} \{ s_1, s_2 \}$.  By standard Sobolev
inequalities, we also have $\nabla_{u^{\epsilon}} u^{\epsilon}$
belongs $L^{\infty}( (s,\infty), H^1)$.  Therefore,
$$||\mathcal{P}^{\a}(\nabla_{u^{\epsilon}} 
u^{\epsilon})||_{\mathcal{V}_{\mu}^1} \leq
||\nabla_{u^{\epsilon}} u^{\epsilon}||_{H^1},$$  implies that for
all $s>0$,  $\mathcal{P}^{\a}(\nabla_{u^{\epsilon}} u^{\epsilon})$
belongs to $L^{\infty}( (s,\infty), \mathcal{V}_{\mu}^1)$, independently of
$\epsilon$. From Lemma \ref{Ulemma} and our $H^2$ estimate
(\ref{h2bd}), we infer that $\mathcal{U}^{\a}(u^{\epsilon})$ is
contained in $L^{\infty}((s,\infty), H^{1+\sigma})$ for
$0<\sigma<1/3$ and consequently
$\mathcal{P}^{\a}(\mathcal{U}^{\a}(u^{\epsilon}))$ belongs
$L^{\infty}((s,\infty), \mathcal{V}_{\mu}^{1+\sigma})$, independently of
$\epsilon$.

Thus using the equation LANS-3 and the results for
$u^{\epsilon}_t$, we have that for any $s>0$,
$$\nu Au^{\epsilon} = - \mathcal{P}^{\a}\left[\nabla_{u^{\epsilon}} u^{\epsilon}
 + \mathcal{U}^{\a}(u^{\epsilon}) -
(1-\a^2 \lap)^{-1}f^{\epsilon} \right] + u^{\epsilon}_t$$ belongs
to $L^{\infty}((s,\infty), \mathcal{V}_{\mu}^1)$, independently of
$\epsilon$.  By the elliptic regularity of the Stokes operator $A$, we
conclude that
$$u^{\epsilon} \in L^{\infty}((s,\infty),
\dot{\mathcal{V}}_{\mu}^3).$$

Define $H = \{u \in L^2 \, | \, \operatorname{div} u = 0 \}$.  
For all $T>0$, (\ref{utbnd}) implies $u_t^{\epsilon} \in
L^2([0,T], H)$ (independently of $\epsilon$) and therefore the
classical compactness theorem (see for instance~\cite{e98,l69})
enables us to conclude that there is a subsequence
$u^{\epsilon'}$ such that for all $T>0$,
\begin{eqnarray*}
u^{\epsilon'} &\rightarrow& u \hspace{.2in} \mathrm{weakly \, in}
\hspace{.2in}
L^2([0,T], D(A)), \\
u^{\epsilon'} &\rightarrow& u \hspace{.2in} \mathrm{strongly \, in}
\hspace{.2in}
L^2([0,T], \mathcal{V}^1_{\mu}), \\
u^{\epsilon'} &\rightarrow& u \hspace{.2in} \mathrm{in} \hspace{.2in}
 C([0,T], H).
\end{eqnarray*}
It is straightforward to verify that $u$ satisfies the weak
formulation associated with $f$.  Indeed, from an identical
argument provided in \cite{fht02}, we may conclude that $u \in
C([0, \infty), \mathcal{V}^1_{\mu})$. Furthermore, we infer that
$u(0)=u_0$. Hence we are lead to the existence of weak solutions
for the LANS-$\a$ equation. The uniqueness and continuous
dependence on initial data of weak solutions can be proved in the
same classical way as was done in~\cite{fht02} on Page 14 for the
periodic case. \qed

\begin{remark}
A standard contraction mapping argument can be used to show the existence
of LANS-$\alpha$ solutions $u$ in $C([0,T],\mathcal{V}^s_{\mu})$ for  $s>1$.
In the absence of forcing or when $f$ is $C^\infty$, the solution $u$ of
the LANS-$\alpha$ equations is
instantly regularized so that $u \in C^{\infty} ((0,\infty) \times\Omega )$.
\end{remark}


\section{Estimating the dimension of the global attractor.} 
\label{sec:attractor}

Theorem \ref{weakth} is sufficient to define the semi-group $S(t)$ by
\[ S(t):u_0 \in \mathcal{V}^1_{\mu} \mapsto u(t) \in \mathcal{V}^1_{\mu}. \]
We now show the uniform compactness property of the operator $S$ and the
existence of an $H^2$ absorbing set.  By Theorem I.1.1 in~\cite{t97}, this
implies the existence of a maximal compact global $H^1$ attractor.
\subsection{Absorbing sets and attractor.}
In this section we prove the existence of an absorbing set in
$\mathcal{V}^2_{\mu}$.  Usually proving the existence of absorbing
sets amounts to proving a priori estimates; for the LANS-$\a$
equations, these estimates are established in Section 2.3. Thus 
estimate
(\ref{h2bd}) implies
$$\displaystyle{\limsup_{t\rightarrow + \infty} c\a^2||u||_{H^2}^2
\leq \rho_2^2 := \nu \lambda_1 K_2 = 2\lambda_1 ||f||^2_{L^2} +
C(||u||_{H^1})}.$$ We conclude that the ball
$B_{\mathcal{V}^2_{\mu}}(0,\rho_2)$ of $\mathcal{V}^2_{\mu}$,
denoted $\mathcal{B}_2$, is an absorbing set in
$\mathcal{V}^2_{\mu}$ for the semi-group $S(t)$.  We choose
${\rho_2}^{\prime}
> \rho_2$ and denote $\mathcal{B}_0$ the ball $B_{\mathcal{V}_{\mu}^1}(0,
{\rho_2}^{\prime})$.  If $\mathcal{B}$ is any bounded set of
$\mathcal{V}^1_{\mu}$, then from
(\ref{h2bd}), $S(t)\mathcal{B} \subset \mathcal{B}_0$ for
$t \geq t_0(\mathcal{B},{\rho_2}^{\prime})$, where
\[ t_0 = s_1 + \frac{1}{\nu \lambda_1} \log \frac{ ||\nabla u(s_1)||^2_{L^2} +
\a^2 ||Au(s_1)||^2_{L^2} }{{{{\rho_2}^{\prime}}}^2 - \rho_2^2}  . \]
At the same time, this result proves the uniform compactness of
$S(t)$; any bounded set $\mathcal{B}$ of $\mathcal{V}^1_{\mu}$  is
included in such a ball, and for $u_0 \in \mathcal{V}^1_{\mu}$ and
$t \geq t_0$, $u(t)$ belongs to $\mathcal{B}_2$ which is bounded
in $H^2$ and relatively compact in $H^1$. Consequently, by Theorem
I.1.1 in~\cite{t97}, the LANS-$\a$ equations have a nonempty
compact, convex, and connected global attractor in $H^1$,
\[ \mathcal{A^{\a}} = \displaystyle{ \cap_{s>0} \left( \cup_{t\geq s}
S(t)\mathcal{B}_2 \right)}. \]
 
\subsection{Dimension of the attractor.}
The dimension of the three-dimensional attractor for the LANS-$\a$
equations follows the arguments of Theorem 6 in ~\cite{fht02}
(stated below), and leads to the same bound as in the periodic
case.  It is important to note that the bound tends to infinity as
$\a \rightarrow 0$.
\begin{theorem}[Theorem 6 of \cite{fht02}]  The Hausdorff and fractal dimensions
of the global attractor of the LANS-$\a$ equations,
$d_H(\mathcal{A^{\a}})$ and $d_F(\mathcal{A^{\a}})$, respectively,
satisfy:
\[ d_H(\mathcal{A^{\a}}) \leq d_F(\mathcal{A^{\a}}) \leq c' \max
\left\{ G^{4/3} \left(\frac{1}{\a^4 \lambda_1^2} \right)^{2/3}, G^{3/2}
\left( \frac{1}{\a^6 \lambda_1^3} \right)^{3/8} \right\}, \]
where $G=\frac{|f|_{L^2}}{\nu^2 \lambda_1^{3/4}}$ is the Grashoff number,
and $c'>0$ is a constant depending only of the shape of $\Omega$.
\end{theorem}
\ni In two dimensions, \cite{cps02} proved that the dimension of the attractor 
can be bounded independently of $\a$ by the Grashoff number.  Hence, in
two-dimensions, we find a similar bound as for the Navier-Stokes
equations.
\begin{theorem}[Theorem 3 of \cite{cps02}]
The Hausdorff and fractal dimensions of the global attractor of the
LANS-$\a$ equations, $d_H(\mathcal{A})$ and $d_F(\mathcal{A})$, respectively,
satisfy:
$$d_H(\mathcal{A}) \leq d_F(\mathcal{A}) \leq C'G,$$
where $\displaystyle{G=\frac{||f||_{L^2}}{\nu^2 \lambda_1}}$ is the Grashoff
number, and $C' >0$ is a constant depending only on the shape of $\Omega$.
\end{theorem}

\newpage
\pagestyle{myheadings}
\chapter{Local well-posedness of the anisotropic LANS-$\alpha$ 
equations on $\mathbb{T}^3$}
\markright{ \rm \normalsize CHAPTER 3. \hspace{.5in} 
Local well-posedness of anisotropic LANS-$\alpha$ on $\mathbb{T}^3$}
\label{chapter:spwp}
\thispagestyle{myheadings}
\section{The anisotropic LANS-$\alpha$ equations.}

The anisotropic Lagrangian averaged Navier-Stokes (LANS-$\alpha$) equations 
are given on $\Omega \times (0,T)$ by
\begin{subequations}\label{aniso}
    \begin{gather}
(1-\alpha^2 C) \left( \partial_t u + \operatorname{div} u \otimes u - \nu 
\mathbb{P} Cu \right) = - \operatorname{grad} p, \label{ueqn} \\
\operatorname{div} u(t,x) = 0, \\
\partial_t F + \nabla F \cdot u = \nabla u \cdot F + [\nabla u \cdot F]^T,
\label{Feqn} \\
u(0,x) = u_0(x), \, F(0,x)= F_0(x) \geq 0,
\end{gather}
\end{subequations}
where $u(t,x)$ denotes the divergence-free mean velocity vector, $p(t,x)$ 
represents the scalar pressure field, $\nu$ is the kinematic viscosity, 
$F(t,x)$ denotes the 
covariance tensor (a $3 \times 3$ matrix), and the notation 
$F_0 \geq 0$ means the initial covariance matrix is assumed positive 
semi-definite.  In \cite{ms03}, the authors define the 
covariance tensor $F$ as the ensemble average of the tensor product of the 
Lagrangian fluctuation vector $\xi^\prime$ given explicitly as 
\[ F(t,x) = \langle \xi^\prime (t,x) \otimes \xi^\prime (t,x) \rangle.\]
Furthermore,  $Cu := \operatorname{div} [\nabla u \cdot F]$ 
and the operator $\mathbb{P}$ is a projection onto divergence-free vector 
fields.  When $F$ is assumed to be a constant multiple of the identity 
matrix, the anisotropic equations reduce to the isotropic equations of Chapter 
\ref{chapter:iso_weak} for only the averaged fluid velocity.  A brief history 
of the isotropic equations can be found in Chapter \ref{chapter:iso_weak}.

The operator $\mathbb{P}$ can be chosen to be either the Leray projector 
$P: L^2 (\Omega) \rightarrow \{v\in L^2(\Omega) \, | \, \operatorname{div} v=0,
\, v \cdot n = 0 \,\, \mbox{on} \,\, \partial \Omega \}$, or the 
generalized Stokes 
projector $\mathcal{P}^{\alpha}_F:= [P(1-\alpha^2 C)]^{-1} P(1-\alpha^2 C)$ 
defined in detail by the following definition.  
\begin{definition}\label{def:gen_stokes}
For any $v$ in the domain of $(1-\alpha^2 C)$ and $F>0$, we set the 
generalized Stokes 
projector $\mathcal{P}_F^{\alpha} (v) =w$, where $w$ is the solution of the 
generalized Stokes problem
\begin{eqnarray*}
(1-\alpha^2 C)w+\operatorname{grad} p &=& (1-\alpha^2 C)v \\
\operatorname{div} w = 0 & & w=0 \, \mathrm{on} \, \partial \Omega;
\end{eqnarray*}
thus 
\[ \mathcal{P}_F^{\alpha}(v) = v - (1-\alpha^2 C)^{-1} \operatorname{grad} p.\]
\end{definition}

\begin{remark}
    The Stokes projector defined in Definition \ref{def:stokes} is a special 
    case of Definition \ref{def:gen_stokes} with the condition $F=\mbox{Id}$.
\end{remark}

The generalized Stokes projector $\mathcal{P}_F^{\alpha}$ maps onto 
the space of $H^1$ divergence-free vector fields that vanish on 
$\partial \Omega$.  Since the Lagrangian fluctuations are necessarily zero 
along $\partial \Omega$, we know $F(t,x) = 0$ on $\partial \Omega$ for all 
time $t\geq 0$.  Unlike the Leray projector which only 
assigns the normal component of $u$ to vanish on the boundary,  the 
generalized Stokes projector provides the correct boundary condition 
for $F$.  Therefore on bounded domains, it is important to use this projector.

Modeling incompressible flow on bounded domains is the purpose of the anisotropic 
model.  The first step towards this goal is proving well-posedness of the 
equations on unbounded domains.  In this chapter, we shall restrict our 
attention to the three-dimensional periodic box $\Omega = 
\mathbb{T}^3$.  We define for $s\geq 0$, 
$$H_{\operatorname{per}}^s = \{ u \in [H^s(\mathbb{T}^3)]^3 \, | \, 
\int_{\mathbb{T}^3} u (x) \, dx = 0 
\}$$ and $$H_{\operatorname{div}}^s = \{ u \in 
H_{\operatorname{per}}^s \, | \, \operatorname{div} u = 0 \}.$$ 
Marsden and Shkoller \cite{ms03} prove the local 
well-posedness of the anisotropic LANS-$\alpha$ equations with the projector 
$\mathbb{P}=P$ for initial velocity fields in the class 
$H_{\operatorname{div}}^s$ and 
positive $F_0 \in [H^s_{\operatorname{per}} (\mathbb{T}^3)]^{3 \times 3}$ for 
$s>3.5$.  They use the classical Galerkin projection to obtain smooth 
finite-dimensional approximations $u_m$ and $F_m$.  After establishing a 
priori bounds independent of $m$, the authors appeal to compactness results to 
extract a subsequence of $u_m$ which converges in 
$C^0([0,T];H_{\operatorname{div}}^s)$ uniformly in $m$.  An important tool in 
obtaining the proper estimates is the control of the regularity of the 
covariance tensor $F$ by the one derivative higher regularity of the 
mean velocity $u$.  This allows the authors to establish estimates 
only dependent on the initial covariance tensor.  

The next step towards applying the anisotropic equations to 
bounded domains is the study of solutions to the equations with the 
generalized Stokes projector viscosity term.   As discussed above, 
the generalized Stokes projector preserved the zero boundary 
condition for the covariance tensor.  Although we known the correct 
projection for bounded domains, the appropriate additional boundary 
condition necessary to solve the forth order equation (\ref{aniso}) is 
unknown.  In this chapter, we show the 
well-posedness of the anisotropic LANS-$\alpha$ equations on 
$\mathbb{T}^{3}$ with the projector 
$\mathbb{P}= \mathcal{P}_F^{\alpha}$ by a method similar to the one used by 
Marsden and Shkoller \cite{ms03}.  However, in obtaining a priori 
$H^s$-type bounds, it is important that we use the differential operator 
$D^{\sigma}P(1-\alpha^2 C)$, $|\sigma| \leq s-2$, rather than the standard 
$D^s$, where $D^s$ denotes all partials derivatives of order $s$.  This 
particular approach is different than what is done in \cite{ms03}.  It allows 
us to estimate the viscosity term without commuting the operator $D^s$ with the 
inverse operator $[P(1-\alpha^2C)]^{-1}$ and leads us to bounds involving the 
Sobolev norms of $P(1-\alpha^2 C)u$ and $u$.   We then appeal to the 
elliptic regularity estimates established in Proposition \ref{ellreg} to write 
all bounds in terms of $P(1-\alpha^2 C)u$.  The same compactness argument 
leads us to the existence of solutions to the anisotropic LANS-$\alpha$
equations (\ref{aniso}).

\section{Local well-posedness with Stokes projector.}

Before proving the existence and uniqueness of classical solutions to the 
anisotropic LANS-$\alpha$ equations stated above, we establish an 
elliptic regularity estimate for the operator $P(1-\alpha^2 C)$ which holds 
assuming the least regularity of $F$.  The proof relies on a standard 
commutator estimate which we include below for completeness.  For the 
remainder of the chapter, $C$ represents an arbitrary positive constant unless, 
by the context, it is clear that $C$ is the operator defined above.

\begin{lemma} 
    For $f, g \in H^{s}(\Omega)$,
\begin{equation}\label{eq:comm}
    ||[D^{s},f]g||_{L^{2}} \leq C \left( ||f||_{H^{s}} 
    ||g||_{L^{\infty}} + ||\nabla f||_{L^{\infty}}||g||_{H^{s-1}} 
    \right)
\end{equation}
\end{lemma}

\begin{proposition}\label{ellreg}
For $F>0$ and $F \in \left[ H^{k-1}_{\operatorname{per}}(\mathbb{T}^3) 
\right]^{3 \times 3}$, suppose $u\in H^k (\mathbb{T}^3)$ solves 
$P(1-\alpha^2 C)u = w$ for $k\geq 3$.  If $w \in H^1 (\mathbb{T}^3)$, then 
\begin{equation}\label{H3_reg}
||u||_{H^3}^2 \leq C (1+C||F||_{H^2}^{16})||w||_{H^1}^2,
\end{equation}
and for $w \in H^{k-2} (\mathbb{T}^3)$, 
\begin{equation}\label{Hk_reg}
||u||_{H^k}^2 \leq C ||w||_{H^{k-2}}^2 + C||F||_{H^{k-1}}^2 ||u||_{H^{k-1}}^2,
\hspace{.2in} k \geq 4.
\end{equation}
Moreover, if $F$ solves (\ref{aniso}) then for some $T>0$,   
\begin{eqnarray}\label{F0_reg}
||u||_{H^3}^2 &\leq& C(t,||F_0||_{H^2})||w||_{H^1}^2 + 
C(||F_0||_{H^2})||w||_{H^1}^4, \label{h3_reg} \\
||u||_{H^4}^2 &\leq&  C(t,||F_0||_{H^2})||w||_{H^2}^2 + 
C(||F_0||_{H^3})||w||_{H^1}^4, \label{h4_reg} \\
||u||_{H^k}^2 &\leq& C ||w||_{H^{k-2}}^2 + 
C(t)||F_0||_{H^{k-1}}^2 ||u||_{H^{k-1}}^2 + 
C||F_0||_{H^{k-1}}^4||u||_{H^{k-1}}^4, 
\end{eqnarray}
for $k\geq 4$ and all $t<T$.
\end{proposition}

\begin{remark} 
    In a classical elliptic regularity statement one 
assumes $u \in H^1$ solves the elliptic second order differential equation  
and attempts to prove $u$ has two derivatives more regularity than the 
inhomogeneous term using difference quotients. Since we are assuming that 
$||F||_{L^\infty}$ is 
not necessarily bounded, this regularity is not guaranteed and we add the 
assumption $u\in H^k$.  Later in the chapter, the proposition is used when 
$u \in C^\infty ([0,T] \times \mathbb{T}^{3})$ and therefore the regularity of 
$u$ is not delicate.
\end{remark}

\begin{remark}
    If $F\in \left [C^{\infty}(\mathbb{T}^{3})\right]^{3 \times 3}$, 
    then equations (\ref{H3_reg}) and (\ref{Hk_reg}) reduce to the 
    standard elliptic regularity result
    \begin{equation}\label{Fsmooth}
        ||u||_{H^{k}}^{2} \leq C(||F||_{W^{k-1,\infty}}) ||w||_{H^{k-1}}^{2}.
    \end{equation}
\end{remark}

\noindent {\bf Proof.} We first write $P(1-\alpha^2 C)u=w$ in the energy form
\begin{equation}\label{energy_eq}
(u,v)_{L^2} + \alpha^2 (F \cdot \nabla u, \nabla v)_{L^2} = (w, v)_{L^2}, 
\hspace{.2in} v \in H_{\operatorname{per}}^k(\mathbb{T}^3) \hspace{.2in} 
k \geq 1.
\end{equation}
We begin by showing the special case of $k=3$, namely inequality 
(\ref{H3_reg}).  When $v=u$, uniform ellipticity of F implies 
$(F \cdot \nabla u, \nabla u)_{L^2} \geq \lambda ||u||_{H^1}^2$ from which we 
find

\begin{equation}\label{H1est}
||u||_{H^1}^2 \leq C ||w||_{H^{-1}}^2.
\end{equation}
Now suppose $v=D^4 u$ and input this into 
inequality (\ref{energy_eq}).  Integrating by parts gives
\begin{equation*}
 \alpha^2 (D^2(F \cdot \nabla u), D^{2} \nabla u)_{L^2} \leq (w, D^4 u)_{L^2}.
\end{equation*}
Expanding the left hand side,
\begin{eqnarray*} 
\alpha^2 (D^2(F \cdot \nabla u), D^2 \nabla u)_{L^2} &=& \alpha^2 (F \cdot D^2 
\nabla u, D^2 \nabla u)_{L^2} + 2 \alpha^2 (\nabla F \cdot D \nabla u, 
D^2 \nabla u)_{L^2} \\ & & +
\alpha^2 (D^2F \cdot \nabla u, D^2 \nabla u)_{L^2} =: I + II + III,
\end{eqnarray*}
where
$I \geq \alpha^2 \lambda ||u||_{H^3}^2$ by uniform ellipticity of F and
\begin{eqnarray*}
II &\leq& C ||\nabla F||_{L^4} ||D^2 u||_{L^4} ||D^3 u ||_{L^2} \leq C 
||F||_{H^2}||D^2 u||_{L^4}||u||_{H^3}, \\
III &\leq& C ||D^2 F||_{L^2} ||\nabla u||_{L^{\infty}}||D^3 u||_{L^2} \leq 
C ||F||_{H^2}||u||_{W^{2,4}}||u||_{H^3} \\ 
&\leq& C ||F||_{H^2} ||D^2 u||_{L^4}||u||_{H^3},
\end{eqnarray*}
by standard Sobolev inequalities.  Combining the estimates above, we get
\begin{eqnarray*}
||u||_{H^3}^2 &\leq&  C (w, D^4u) + C ||F||_{H^2}||D^2 
u||_{L^4}||u||_{H^3} \\
&\leq&  C ||w||_{H^1}||D^4u||_{H^{-1}} + C ||F||_{H^2}||u||_{H^3}^{7/8} 
||u||_{H^1}^{1/8} ||u||_{H^3} \\
&\leq& \varepsilon ||u||_{H^3}^2 + C ||w||_{H^1}^2 + 
C||F||_{H^2}^{16} ||u||_{H^1}^2,
\end{eqnarray*}
where we have used the Gagliardo-Nirenberg inequality (\ref{gn}) for the 
second inequality and Young's inequality (\ref{eq:youngs}) for the 
last.  Using estimate (\ref{H1est}) 
and choosing $\varepsilon >0$ sufficiently small, we can absorb the 
$||u||_{H^3}^2$-term into the left side, obtaining
\[ ||u||_{H^3}^2 \leq C (1+C||F||_{H^2}^{16})||w||_{H^1}^2. \]

For $k\geq 4$, we substitute $v=(-1)^{k-1}D^{2k-2} u$ into (\ref{energy_eq}). 
After integrating by parts,
\[ \alpha^2 (D^{k-1} (F \cdot \nabla u), D^{k-1} \nabla u)_{L^2} \leq 
(D^{k-1} w, D^{k-1} u)_{L^2}. \]
Expanding the left hand side we see 
\begin{eqnarray*}
\alpha^2 (D^{k-1} (F \cdot \nabla u), D^{k-1} \nabla u)_{L^2} &=& \alpha^2 (F 
\cdot D^{k-1} \nabla u, D^{k-1} \nabla u)_{L^2} \\ 
& & + \alpha^2 ([D^{k-1}, F] 
\nabla u, D^{k-1} \nabla u)_{L^2} =: I + II, 
\end{eqnarray*}
where $I \geq \alpha^2 \lambda ||u||_{H^k}^2$. Unlike the case when $k=3$, 
$k \geq 4 > \frac{7}{2}$ implies $H^{k-1} \subset W^{1,\infty}$ and therefore 
by a standard commutator estimate (\ref{eq:comm}),
\[ II \leq C (||F||_{H^{k-1}}||\nabla u||_{L^{\infty}} + 
||\nabla F||_{L^{\infty}} ||\nabla u||_{H^{k-2}})||u||_{H^k} \leq C 
||F||_{H^{k-1}}||u||_{H^{k-1}}||u||_{H^k}. \]
Hence 
\begin{eqnarray*}
||u||_{H^k}^2 &\leq& C ||D^{k-1}w||_{H^{-1}}||D^{k-1}u||_{H^1} + C 
||F||_{H^{k-1}}||u||_{H^{k-1}}||u||_{H^k}, \\
&\leq& \varepsilon ||u||_{H^k}^2 + C ||w||_{H^{k-2}}^2 + C ||F||_{H^{k-1}}^2 
||u||_{{H^{k-1}}}^2.
\end{eqnarray*}
By choosing $\varepsilon >0$ sufficiently small we achieve inequality 
(\ref{Hk_reg}).  

To show inequality (\ref{h3_reg}), recall from \cite{ms03}, equation 
(33), that if $F$ solves (\ref{Feqn}) for a given $u$, we may estimate 
$||F||_{H^s}$ by one higher order derivative of $u$ in the following 
way, 
\begin{equation}\label{F_bound}
||F||_{H^s} \leq C ||F_{0}||_{H^{s}} \exp \int_{0}^{t} 
||u(s)||_{H^{s+1}} \, ds \leq C ||F_0||_{H^s} \left( 1 + t ||u||_{H^{s+1}} + O(t^2) \right).
\end{equation}
For $k=3$, it follows that inequality (\ref{H3_reg}) can be written as
\begin{eqnarray*}
||u||_{H^3}^2  &\leq& C||w||_{H^1}^2 + C||F_0||_{H^2}^{16} \left( 1+ 
16t||u||_{H^3} + O(t^2) \right) ||w||_{H_1}^2, \\ 
&\leq& C||w||_{H^1}^2 + C||F_0||_{H^2}^{16} \left(1+ O(t^2)\right) 
||w||_{H^1}^2+ tC||F_0||_{H^2}^{16}||w||_{H^1}^2||u||_{H^3}, \\
&\leq&  C||w||_{H^1}^2 + C(t) ||F_0||_{H^2}^{16} ||w||_{H^1}^2+ 
C||F_0||_{H^2}^{32}||w||_{H^1}^4 + \varepsilon t^2||u||_{H^3}^2,
\end{eqnarray*}
For $t<T$, we choose $\varepsilon >0$ sufficiently small to conclude 
inequality (\ref{F0_reg}).  The estimate for $k \geq 4 $ follows from the same 
procedure.
\qed

\begin{theorem}\label{stokes_wp} For $s>7/2$, and $u_0 \in H^s_{\operatorname{div}} 
(\mathbb{T}^3)$, $F_0 \in [H^s_{\operatorname{per}} 
(\mathbb{T}^3)]^{3 \times 3}$, with $F_0>0$, there exists a unique solution 
$(u, F)$ with $u \in C^0([0,T]; H^s_{\operatorname{div}}) \cap L^2 (0,T; 
H^{s+1}_{\operatorname{div}})$ and $F \in C^0([0,T]; 
[H^s_{\operatorname{per}}]^{3 \times 3})$ to the anisotropic LANS-$\alpha$ 
equations, where $T$ depends on the initial data.
\end{theorem}

\noindent {\bf Proof.}  The evolution equation for the mean velocity 
(\ref{ueqn}) may be written as 
\begin{eqnarray}\label{aniso2}
\partial_t u + \mathcal{P}^{\alpha}_F \operatorname{div}(u \otimes u) &=& \nu 
\mathcal{P}^{\alpha}_F Cu \\
u(0,x) &=& u_0(x). \nonumber
\end{eqnarray}

\noindent {\bf Approximate solutions.}  Let $v_i = v_{i}(\cdot, x)$ 
$(i=1, 2, \dots)$ denote the smooth periodic orthogonal basis of 
$H^1_{\operatorname{per}}$ given by the eigenfunctions of the Stokes 
operator.  Define 
$P_m:H^1_{\operatorname{per}} \rightarrow V^{m}:=\operatorname{span} \{ 
v_{1}, \dots, v_{m} \}$ and $u_{m} = P_{m} u$. 


Consider the Galerkin projection of (\ref{aniso}) given by
\begin{equation}\label{galerkin}
    \begin{split}
\partial_t u_m + P_m \mathcal{P}^{\alpha}_{F_m} \operatorname{div}(u_m 
\otimes u_m) &= \nu P_m \mathcal{P}^{\alpha}_{F_m} C_{m} u_m, \\
u_m(0,x) &= P_m u_0,
\end{split}
\end{equation}
with $C_{m}u_m := \operatorname{div} (F_{m} \cdot \nabla u_{m})$ 
where $F_{m}$ is a solution to
\begin{equation}\label{eq:Fm}
    \begin{split}
\partial_t F_{m}  + \nabla F_{m} \cdot u_{m} &= \nabla u_{m} \cdot F_{m} + 
[\nabla u_{m} \cdot F_{m}]^T, \\
F_{m}(0,x) &= F_0(t,x) > 0.
\end{split}
\end{equation}

For each $m \in \mathbb{N}$, there is a smooth solution; in order to 
pass to the limit as $m\rightarrow \infty$  to produce a solution 
to (\ref{aniso}), it suffices to obtain a priori estimates independent 
of $m$.

\noindent {\bf $H^{s}$ estimate.}  Since each term on the right-hand side of 
equation (\ref{galerkin}) involves the generalized Stokes projector 
$\mathcal{P}_{F_{m}}^{\alpha} := \left[ P(1-\alpha^2 C_{m})\right]^{-1} 
P(1-\alpha^2 C_{m})$, we use the operator $P(1-\alpha^2 C_{m})$ before the 
standard differential operator in the following $H^{s}$-energy estimate.
By letting $\sigma$ denote a multi-index with $|\sigma| \leq s-2$, we apply 
$D^{\sigma}P(1-\alpha^2 C_{m})$ to both sides of equation (\ref{galerkin}) 
and take the $L^2$-inner product with $D^{\sigma} P (1-\alpha^2 
C_{m})u_m$.  The left-hand side reduces to
\begin{align*}
(D^{\sigma}P(1-\alpha^2 C_{m})\partial_t u_m &, D^{\sigma} w_m)_{L^2} = 
(D^{\sigma} \partial_t w_m, D^{\sigma} w_m)_{L^2} - \alpha^2 (D^{\sigma} 
P[C_{m}, 
\partial_t] u_m, D^{\sigma} w_m)_{L^2} \\
&= (\partial_t D^{\sigma} w_m, D^{\sigma} w_m)_{L^2} +\alpha^2 
(D^{\sigma} P \operatorname{div}(\partial_t F_{m} \cdot \nabla u_m), D^{\sigma} 
w_m)_{L^2} 
\\ &= \frac{1}{2} \frac{d}{dt} ||D^{\sigma} w_m||_{L^2}^2 - \alpha^2 
(D^{\sigma}P(\partial_t F_{m} \cdot \nabla u_m), \nabla D^{\sigma} 
w_m)_{L^2},
\end{align*}
where we define $w_m : =P(1-\alpha^2 C_{m})u_m$.
In addition
\begin{eqnarray}\label{w'}
\frac{d}{dt} ||D^{\sigma} w_m||_{L^2}^2 &=& 2\alpha^2 
(D^{\sigma}P(\partial_t F_{m} \cdot \nabla u_m), \nabla D^{\sigma} w_m)_{L^2} 
\nonumber \\ 
& & - 2 (D^{\sigma} P (1-\alpha^2 C_{m}) \operatorname{div}(u_m \otimes u_m), 
D^{\sigma} w_m)_{L^2} \nonumber \\ 
& & + 2 \nu (D^{\sigma} P(1-\alpha^2 C_{m})C_{m}u_m, 
D^{\sigma} w_m)_{L^2} \nonumber \\
&=&  I + II + III,
\end{eqnarray}
where we have used the definition of $\mathcal{P}^{\alpha}_{F_{m}} 
:= [P(1-\alpha^2 C_{m})]^{-1} P(1-\alpha^2 C_{m})$ for the $II$ and $III$ terms. 

\noindent {\bf Estimate for I.} By the continuity of the Leray projector $P$ 
and the fact that  $H^r$ is an multiplicative algebra for $r>\frac{3}{2}$, we
use equation (\ref{eq:Fm}) to write the first term as 
\begin{eqnarray*}
|I| &\leq& C(\alpha^2) ||P(\partial_t F_{m} \cdot \nabla u_m)||_{H^{s-2}} 
||w_m||_{H^{s-1}} \leq C(\alpha^2) ||\partial_t F_{m}||_{H^{s-2}} 
||u_m||_{H^{s-1}} ||w_m||_{H^{s-1}} \\
&\leq& C(\alpha^2) (||\nabla F_{m} \cdot u_m ||_{H^{s-2}} + ||\nabla u_m \cdot 
F_{m}||_{H^{s-2}} )||u_m||_{H^{s-1}} ||w_m||_{H^{s-1}} \\
&\leq&  
C(\alpha^2) ||F_{m}||_{H^{s-1}} ||u_m||_{H^{s-1}}^2 ||w_m||_{H^{s-1}} \\
&\leq&  \varepsilon_1 ||w_m||_{H^{s-1}}^2 + C(\varepsilon_1, \alpha^4) 
||F_{m}||_{H^{s-1}}^2 ||u_m||_{H^{s-1}}^4,
\end{eqnarray*}
where we have used Young's inequality for the last estimate.

\noindent {\bf Estimate for II.} Since $\operatorname{div} D^{\sigma}w_m = 
0$, 
\begin{eqnarray*}
|II| &\leq& C|(D^{s-2} (1-\alpha^2 C_{m})\operatorname{div}(u_m \otimes u_m), 
D^{s -2}w_m)| \\
&\leq& C ||\operatorname{div} (u_m \otimes 
u_m)||_{H^{s-2}}||w_m||_{H^{s-2}} \\ 
& & +
C(\alpha^2) |(D^{s-2} (F_{m} \cdot \nabla \operatorname{div} (u_m \otimes u_m)), 
\nabla D^{s-2} w_{m})| \\
&\leq & C||\nabla u_m \cdot u_m ||_{H^{s-2}} ||w_m||_{H^{s-2}} + C(\alpha^2)
||F_{m} \cdot \nabla \operatorname{div} (u_m \otimes u_m))||_{H^{s-2}} 
||w_m||_{H^{s-1}}. 
\end{eqnarray*}
By Young's inequality we conclude 
\[ |II| \leq \varepsilon_2 ||w_m||_{H^{s-1}}^2 + C ||u_m||_{H^{s-2}} 
||u_m||_{H^{s-1}} ||w_m||_{H^{s-2}}+ C(\alpha^4) ||F_{m}||_{H^{s-2}}^2 
||u_m||_{H^{s-1}}^2 ||u_m||_{H^{s}}^2. \]

\noindent {\bf Estimate for III.} 
Using the coercivity of $F_{m}$, the estimate for expression III 
includes a negative coefficient in front of the 
$||w_{m}||_{H^{s-1}}$-term.  The negative coefficient allows us to control the 
$||w_{m}||_{H^{s-1}}$-terms from estimates for expressions I and II.

\begin{eqnarray}\label{III}
III &=& 2 \nu (D^{\sigma}P(1-\alpha^2 C_{m}) C_{m}u_m, D^{\sigma}w_m)_{L^2} = 
2 \nu (D^{\sigma} P C_{m} (1-\alpha^2 C_{m})u_m, D^{\sigma} w_m)_{L^2} \nonumber \\
&=& 2 \nu (D^{\sigma} C_{m}P(1-\alpha^2 C_{m})u_m, D^{\sigma}w_m)_{L^2} + 2 \nu 
(D^{\sigma} [P,C_{m}](1-\alpha^2 C_{m}) u_m, D^{\sigma } w_m)_{L^2} \nonumber \\
&=& -2 \nu (D^{\sigma}(F_{m} \cdot \nabla w_m), \nabla D^{\sigma}w_m)_{L^2} 
+2\nu (D^{\sigma} [P,C_{m}]u_m, D^{\sigma} w_m)_{L^2} \nonumber \\
& & - 2 \alpha^2 \nu (D^{\sigma}[P,C_{m}]C_{m}u_m, D^{\sigma} w)_{L^2}.
\end{eqnarray}
Since $\operatorname{div} u_m = 0$, the definition of the Leray projector 
$P$ implies $[P,C_{m}]u_m=PC_{m}u_m-C_{m}Pu_m=- \operatorname{grad} p_1$ 
where $\lap p_1 =  \operatorname{div} C_{m}u_m$.  Therefore the second term on 
the right hand side of equation (\ref{III}) vanishes.  Also  
\begin{eqnarray*} 
[P,C_{m}]C_{m}u_m &=& PC_{m}^2u_m-C_{m}PC_{m}u_m= C_{m}^2u_m - 
\operatorname{grad} p_2 - C_{m}^2u_m + C_{m} 
\operatorname{grad} p_1 \\
&=& - \operatorname{grad} p_2 + \operatorname{grad} C_{m}p_1 + 
[C_{m}, \nabla]p_1 =: 
\operatorname{grad} P - \operatorname{div}(\nabla F_{m} \cdot \nabla p_1),
\end{eqnarray*}
where $\lap p_2 = \operatorname{div} C_{m}^2 u_m$.  
Using Leibniz' formula on the first term of (\ref{III}) we find
\begin{eqnarray}\label{IIIest}
III &=& -2 \nu (F_{m} \cdot \nabla D^{\sigma} w_m, \nabla D^{\sigma} 
w_m)_{L^2} +C\displaystyle{\sum_{0<|\beta| \leq |\sigma|}}(D^{\beta}F_{m} 
D^{\sigma- \beta} \nabla w_m, \nabla D^{\sigma} w_m)_{L^2} \nonumber \\
& & - 2 \alpha^2 \nu (D^{\sigma}(\nabla F_{m} \cdot \nabla p_1), 
\nabla D^{\sigma} w_m)_{L^2}. 
\end{eqnarray}
Rather than introducing the commutator in equation (\ref{IIIest}), we use the 
Leibniz' formula to obtain the correct power on the term $||w_m||_{H^3}$ for 
the critical case of $\sigma=4$. With a standard commutator estimate the best 
we can obtain is a bound by $||F_{m}||_{H^3}||w_m||_{H^3}^2$, which doesn't allow 
the correct isolation of the $||w_m||_{H^3}$-term.    

By the uniform ellipticity of $F_{m}$, the first term of (\ref{IIIest})
\[  -2 \nu (F_{m} \cdot \nabla D^{\sigma} w_m, \nabla D^{\sigma}w_m)_{L^2}  
\leq -\lambda ||w_m||_{H^{s-1}}^2. \]
The second term of inequality (\ref{IIIest})
\begin{align*}
\displaystyle{\sum_{0<|\beta| \leq |\sigma|}} 
(D^{\beta}F_{m} D^{\sigma-\beta} 
\nabla w_m &, \nabla D^{\sigma} w_m)_{L^2} \leq C 
\displaystyle{\sum_{0<|\beta |\leq |\sigma|}}||D^{\beta} F_{m}||_{L^4} 
||D^{\sigma -\beta +1}w_m||_{L^4} ||\nabla D^{\sigma}w_m||_{L^2} \\
& \leq C ||D^{s-2}F_{m}||_{H^1} ||D^{s - 2} w_m ||_{L^4} 
||w_m||_{H^{s-1}} \\
& \leq C ||F_{m}||_{H^{s-1}} ||D^{s - 1}w_m||_{L^2}^{\frac{3}{4}} 
||D^{s-2} w_m||_{L^2}^{\frac{1}{4}} ||w_m||_{H^{s-1}} \\ 
&\leq C ||F_{m}||_{H^{s-1}} ||w_m||_{H^{s-2}}^{\frac{1}{4}} 
||w_m||_{H^{s-1}}^{\frac{7}{4}}, 
\end{align*}
where we have used the Gagliardo-Nirenberg inequality for the third 
inequality.  Finally the last term of (\ref{IIIest}) is 
\begin{eqnarray*}
- 2 \alpha^2 \nu (D^{\sigma}(\nabla F_{m} \cdot \nabla p_1), 
\nabla D^{\sigma} w_m)_{L^2} &\leq& C(\alpha^2) ||F_{m}||_{H^{s-1}} 
||\nabla p_1||_{H^{s-2}} ||w_m||_{H^{s-1}}  \\
&\leq& C(\alpha^2) ||F_{m}||_{H^{s-1}}||C_{m}u_m||_{H^{s-2}}||w_m||_{H^{s-1}} \\
&\leq& C(\alpha^2) ||F_{m}||_{H^{s-1}}^2 ||u_m||_{H^s} ||w_m||_{H^{s-1}}
\end{eqnarray*}
Hence by Young's inequality we combine the above estimates to write 
inequality  (\ref{IIIest}) as
\begin{equation*}
|III| \leq (-\lambda + \varepsilon_3)||w_m||_{H^{s-1}}^2 + C 
||F_{m}||_{H^{s-1}}^8 
||w_m||_{H^{s-2}}^2 + C (\alpha^2) ||F_{m}||_{H^{s-1}}^4 ||u_m||_{H^{s}}^2
\end{equation*}

\noindent Finally, the estimates for $I$, $II$, and $III$ allow us to conclude
\begin{eqnarray}\label{w'est}
\frac{d}{dt} ||w_m||_{H^{s-2}}^2 &\leq& (-\lambda+ \varepsilon) 
||w_m||_{H^{s-1}}^2 + C||F_{m}||_{H^{s-1}}^8 ||w_m||_{H^{s-2}}^2 \nonumber \\
& & +C ||u_m||_{H^{s-2}} ||u_m||_{H^{s-1}}||w_m||_{H^{s-2}} + C 
||F_{m}||_{H^{s-2}}^2 ||u_m||_{H^{s-1}}^2||u_m||_{H^s}^2 \nonumber \\
& & +C ||F_{m}||_{H^{s-1}}^4||u_{m}||_{H^s}^2 + 
C||F_{m}||_{H^{s-1}}^2 ||u_m||_{H^{s-1}}^4.
\end{eqnarray}

Using inequality (\ref{F_bound}) and the 
elliptic regularity estimates of Proposition \ref{ellreg}, we can bound the 
right hand side of (\ref{w'est}) by 
only $t$, $||F_0||_{H^{k}}$, and $||w_m||_{H^k}$ for $k\leq s-1$.   
This is precisely the bound we need to show the existence of classical 
solutions using Proposition \ref{ellreg}. As an example of this 
computation, we include the case when $s=4$.  

\noindent {\bf Example $s = 4$.}  We now use Proposition 
\ref{ellreg} and inequality (\ref{F_bound}) to write the estimate above for 
$s=4$, in terms of $w_m$ and $F_0$ alone.  We begin with the second term,
\begin{eqnarray*}
||F_{m}||_{H^3}^8 ||w_m||_{H^2}^2 &\leq& C ||F_0||_{H^3}^8 (1+8t||u_m||_{H^4} + 
O(t^2))||w_m||_{H^2}^2 \\
&\leq& C ||F_0||_{H^3}^8 (1+O(t^2))||w_m||_{H^2}^2 + tC||F_0||_{H^3}^8 
||u_m||_{H^4} ||w_m||_{H^2}^2 \\
&\leq& C(t,||F_0||_{H^3}) ||w_m||_{H^2}^2 + C(||F_0||_{H^3})||w_m||_{H^2}^4 + 
C(t, ||F_0||_{H^3}) ||u_m||_{H^4}^2 \\
&\leq&  C(t,||F_0||_{H^3}) ||w_m||_{H^2}^2 + C(t,||F_0||_{H^3})||w_m||_{H^2}^4,
\end{eqnarray*} 
where we have used Young's inequality for the third inequality and inequality 
(\ref{h4_reg}) for the last.  The third term in inequality (\ref{w'est}) can 
be estimated by a repeated use of Young's inequality,
\begin{eqnarray*}
||u_m||_{H^2} ||u_m||_{H^3} ||w_m||_{H^2} &\leq& C ||w_m||_{H^2}^2 + C 
||u_m||_{H^3}^4 \\
&\leq& C ||w_m||_{H^2}^2 + C (t, ||F_0||_{H^3}) ||w_m||_{H_2}^4 \\
& & + C(t,||F_0||_{H^3}) ||w_m||_{H^1}^4||w_m||_{H^2}^2 + 
C(||F_0||_{H^3})||w_m||_{H^1}^8 \\ 
&\leq&  C ||w_m||_{H^2}^2 + C (t, ||F_0||_{H^3}) ||w_m||_{H^2}^4 \\
& & + C(t,||F_0||_{H^3})||w_m||_{H^1}^8.
\end{eqnarray*}
The fourth term of inequality (\ref{w'est}) can be estimated by
\begin{eqnarray*}
||F_{m}||_{H^2}^2 ||u_m||_{H^3}^2 ||u_m||_{H^4}^2 &\leq& C(t)||F_0||_{H^2}^2 
||u_m||_{H^4}^4 + tC||F_0||_{H^2}^2 ||u_m||_{H^3}^3 ||u_m||_{H^4}^2 \\
&\leq& C(t,||F_0||_{H^3})||u_m||_{H^4}^4 + C(t,||F_0||_{H^2})||u_m||_{H^3}^6 \\
&\leq& C(t,||F_0||_{H^2})||w_m||_{H^1}^4 + C(t,||F||_{H^2}) 
||w_m||_{H^1}^8 \\ & & +  
C(t,||F_0||_{H^3}) ||w_m||_{H^2}^4,
\end{eqnarray*}
by inequalities (\ref{h3_reg}), (\ref{h4_reg}), and Young's inequality.  The 
second to last term of inequality (\ref{w'est})
\begin{eqnarray*}
||F_{m}||_{H^3}^4 ||u_m||_{H^4}^2 &\leq& C(t) ||F_0||_{H^3}^4 ||u_m||_{H^4}^2 + t 
C||F_0||_{H^3}^4 ||u_m||_{H^4}^3 \\
&\leq& C(t,||F_0||_{H^3})||u_m||_{H^4}^2+C(t,||F_0||_{H^3})||u_m||_{H^4}^4 \\ 
&\leq& C(t,||F_0||_{H^3})||w_m||_{H^1}^4 + C(t,||F_0||_{H^3})||w_m||_{H^1}^8 \\
& & + C(t,||F_0||_{H^3})||w_m||_{H^2}^2 + C(t,||F_0||_{H^3})||w_m||_{H^2}^4,
\end{eqnarray*}
where we have used Young's inequality to write 
$||u_m||_{H^4}^3 \leq C||u_m||_{H^{4}}^2 + C||u_m||_{H^4}^4.$
The final term of inequality (\ref{w'est}) can be estimated by inequality 
(\ref{H3_reg}) as 
\begin{eqnarray*}
||F_{m}||_{H^3}^2||u_m||_{H^3}^4 &\leq& C ||F_{m}||_{H^3}^2 (1+C 
||F_{m}||_{H^2}^{16})^2 ||w_m||_{H^1}^4 \\
&\leq& C \left[ ||F_{m}||_{H^3}^2 + ||F_{m}||_{H^3}^{18} + 
||F_{m}||_{H^3}^{34} \right ] ||w_m||_{H^1}^4 \\
&\leq& C(t,||F_0||_{H^3})||w_m||_{H^1}^4 + t C(||F_0||_{H^3})||u_m||_{H^4} 
||w_m||_{H^1}^4 \\
&\leq& C(t,||F_0||_{H^3})||w_m||_{H^1}^4 + C(t,||F_0||_{H^3}) ||w_m||_{H^1}^8 
\\ & & + C(t,||F_0||_{H^3})||u_m||_{H^4}^2 \\
&\leq& C(t,||F_0||_{H^3})||w_m||_{H^1}^4 + C(t,||F_0||_{H^3}) ||w_m||_{H^1}^8 
\\
& & + C(t,||F_0||_{H^3})||w_m||_{H^2}^2.
\end{eqnarray*} 
Finally, combining the estimates above we achieve the a priori estimate
\begin{equation}\label{w_H2_final}
\frac{d}{dt} ||w_m||_{H^2}^2 \leq (-\lambda + \varepsilon)||w_m||_{H^3}^2 + C 
||w_m||_{H^2}^2 + C||w_m||_{H^2}^4 + C||w_m||_{H^1}^4 + C||w_m||_{H^1}^8,
\end{equation}
where $C=C(t, \alpha, ||F_0||_{H^3})>0$.    

\noindent {\bf Convergence to a strong solution to equation (\ref{aniso}).}

From inequality (\ref{w_H2_final}) and equation (\ref{galerkin}), we may 
conclude that for some $T>0$, $w_m$ is bounded in $L^{\infty}(0,T; 
H_{\operatorname{div}}^2) \cap W^{1,\infty} (0,T; L_{\operatorname{div}}^2) 
\cap L^2 (0,T;H_{\operatorname{div}}^3)$ 
uniformly in $m$, and hence by inequality (\ref{h4_reg}), 
\[ u_m \, \, \mbox{is bounded in} \, \, L^{\infty}(0,T; 
H_{\operatorname{div}}^4) \cap W^{1,\infty} (0,T; H_{\operatorname{div}}^2) 
\cap L^2 (0,T;H_{\operatorname{div}}^{5}),\] uniformly in $m$.  By the weak 
compactness theorem, there exists a subsequence $u_{m_k}$ such that for all 
$T>0$,
\[ u_{m_k} \rightharpoonup u\, \, \mbox{in} \, \, L^{\infty}(0,T; 
H_{\operatorname{div}}^4) \cap W^{1,\infty} (0,T; H_{\operatorname{div}}^2) 
\cap L^2 (0,T;H_{\operatorname{div}}^5). \]
Since $W^{1,\infty} (0,T;H_{\operatorname{div}}^2) \subset C^0([0,T]; 
H_{\operatorname{div}}^2)$, Arzela-Ascoli compactness criterion implies 
\[ u_{m_k} \rightarrow u \, \, \mbox{in} \, \, C^0([0,T]; 
H_{\operatorname{div}}^2).\]
Furthermore by interpolation, $u_m \in W^{0,\infty} 
(0,T;H_{\operatorname{div}}^4) \cap W^{1,\infty} (0,T; 
H_{\operatorname{div}}^2) $ implies
\[ u_{m_k} \rightarrow u \, \, \mbox{in} \, \, C^\delta ([0,T]; 
H_{\operatorname{div}}^{4-\delta}), \]
and therefore by Sobolev's embedding theorem (\ref{eq:sobolev}),
\[  u_{m_k} \rightarrow u \,\, \mbox{in} \,\, C^0 ([0,T]; 
C_{\operatorname{div}}^2) \]
for $\delta$ taken sufficiently small.  Thus all the terms on the right-hand 
side of equation (\ref{galerkin}) converge strongly and $\partial_t u_m$ 
converges weakly.  To see that $u \in C^0([0,T]; H_{\operatorname{div}}^4)$, 
it suffices to show that $||u(t,\cdot)||_{H^4}$ is continuous on $[0,T]$, but 
this follows from the inequalities (\ref{h4_reg}) and (\ref{w_H2_final}).

The uniqueness of classical solutions to (\ref{aniso}) follows the same 
arguments made in \cite{ms03}. The only deviation being the evolution 
of the difference of two solutions of (\ref{aniso2}) is given as
\begin{eqnarray*}
\partial_t (u_1-u_2) &+& (\mathcal{P}_1^{\alpha} - \mathcal{P}_2^{\alpha}) 
\nabla_1 u_1 \cdot u_1  + \mathcal{P}_2^{\alpha} \nabla u_1 \cdot (u_1 - u_2) 
+ \mathcal{P}_2^{\alpha} \nabla (u_1-u_2) \cdot u_2  \\
&=& \nu \mathcal{P}_1^{\alpha} \operatorname{div} ((F_1-F_2)\cdot \nabla u_1) 
+ \nu \mathcal{P}_1^{\alpha} \operatorname{div} (F_2 \cdot \nabla u_1) + \nu 
\mathcal{P}_2^{\alpha} \operatorname{div} (F_2 \cdot \nabla u_2) \\
&=&  \nu \mathcal{P}_1^{\alpha} \operatorname{div} ((F_1-F_2)\cdot \nabla u_1) 
+ \nu \mathcal{P}_1^{\alpha} \operatorname{div} (F_2 \cdot \nabla 
(u_1-u_2)) \\ & & + 
\nu (\mathcal{P}_1^{\alpha} -\mathcal{P}_2^{\alpha}) \operatorname{div} 
(F_2 \cdot \nabla u_1), 
\end{eqnarray*}
where $F_i$ solves equation (\ref{Feqn}) with $u=u_i$ for $i=1,2$. That is, if
\[ y(t):=||u_1(t)-u_2(t)||_{H^1}^2 + ||F_1(t)-F_2(t)||_{H^1}^2 \]
we use the fact that $u\in L^\infty(0,T;H_{\operatorname{div}}^s) \cap L^2 
(0,T;H_{\operatorname{div}}^{s+1})$ and computations similar to that used for 
existence to obtain the differential inequality
\[ \frac{d}{dt} y(t) \leq C(t)y(t), \hspace{.2in} C(t)=C(||u_i||_{H^s}, 
||F_i||_{H^s}), \hspace{.2in} i=1,2, \]
from which uniqueness follows.
\qed

\begin{remark}
    As $\alpha \rightarrow 0$, the anisotropic LANS-$\alpha$ 
    equations should reduce to the Navier-Stokes equations.    Once 
    it is shown that $F \rightarrow \mbox{Id}$ as $\alpha 
    \rightarrow 0$, the anisotropic equations reduce to the isotropic 
    equations.  Marsden and Shkoller \cite{ms01} have proven for 
    $s\geq 3$, solutions to the isotropic LANS-$\alpha$ equations converge 
    in $H^{s}$ for short time on intervals which are governed by the 
    existence theory for the Navier-Stokes equations.
\end{remark}

\newpage
\pagestyle{myheadings}
\chapter{Numerical solutions to the anisotropic LANS-$\alpha$ 
equations}
\markright{ \rm \normalsize CHAPTER 4. \hspace{.5in} 
Numerical solutions to the anisotropic LANS-$\alpha$ equations}
\label{chapter:numerical}
\thispagestyle{myheadings}
\section{Introduction.}

Examining the behavior of fluid flow in elementary domains is more than a 
common numerical test for a new fluid model, it is also a means of 
highlighting interesting phenomena inherent in the model.  
The direct numerical simulation (DNS) of turbulent 
flow at small to moderate Reynolds number has become an important 
computational tool in understanding large scale turbulence 
motion.  Unfortunately, DNS are still computationally expensive in 
turbulent regimes.  Alternate 
approaches to brute-force DNS are the Reynolds averaged 
Navier-Stokes (RANS) simulations and large eddy simulations (LES).  The 
RANS equations (\ref{eq:RANS}) are briefly discussed in Chapter 
\ref{chapter:introduction}.
In LES, a spatial averaging operator (filter) is applied to the 
Navier-Stokes equations to obtain a new set of equations for the 
averaged (filtered) variables.  Due to the computational limit of DNS 
for large Reynolds number, LES have become one of the 
standard methods in solving for fluid flow.

The behavior of small spatial scales in 
turbulent flow is often characterized by statistical isotropy 
and homogeneity away from the boundary of the fluid container.  
Therefore, the 
isotropic LANS-$\alpha$ equations 
appear to be an appropriate model in isotropic regimes and have 
been studied recently from the numerical point of view.  As a proposed model 
for large scale turbulence, Mohseni \emph{et. al.} \cite{msk03} 
compare the isotropic LANS-$\alpha$ equations to known results for 
DNS and LES methods.  The authors demonstrate the utility of the isotropic 
~LANS-$\alpha$ as a sub-grid stress model for three-dimensional isotropic 
forced and decaying turbulence.  They perform two sets of 
forced isotropic turbulence simulations and compare their results to 
DNS and LES results where appropriate.   In the LANS-$\alpha$ 
simulations, Mohseni \emph{et. al.} conclude that 
selecting an appropriate $\alpha$ is a delicate 
comprise between the accuracy of the model at the the large scales and 
the minimum resolution requirements.  In particular, the 
LANS-$\alpha$ equations accurately mimic the behavior of the 
Navier-Stokes equations at large spatial scales as long as a minimum 
resolution is observed.  The accuracy improves as $\alpha \rightarrow 
0$, 
but the computation requires higher resolution.  The higher resolution 
comes with a price: an unresolved LANS-$\alpha$ computation could result in 
the contamination of large scales and a loss of accuracy.  

In 1998, Chen \emph{et. al.} \cite{cfh98} studied the mean velocity of 
turbulent channel and pipe flows.  They proposed using stationary 
solutions to isotropic LANS-$\alpha$ equations as a closure
approximation for the Reynolds-averaged equations.  Since in the 
near-wall region the fluctuations are highly anisotropic, their 
results were only in good agreement with experimental data away from the 
viscous boundary layer.  

By restricting the anisotropic LANS-$\alpha$ equations to the 
channel, Coutand and Shkoller \cite{cs03a} propose a turbulent channel theory 
that models the large scale fluid motion throughout the entire domain.  Unlike 
the isotropic equations, the solutions to the anisotropic equations consist 
of both a mean velocity field $u$ and a $3 \times 3$ covariance tensor 
$F$, defined, in detail, in the next section.  An important property of $F$ 
is its degeneracy to zero at the wall.   Coutand 
and Shkoller \cite{cs03a} show that near the 
wall, the degeneracy scales like $d \sqrt{ | \log d |}$, where $d$ is the 
normalized distance function to the wall.  The authors 
compensate for the degeneracy in the boundary layer by working in 
weighted Sobolev spaces.  In this functional framework, they prove the 
global-in-time existence and uniqueness of weak solutions to the 
anisotropic LANS-$\alpha$ equations.  In particular, they restrict the 
anisotropic equations to the channel and make the 
assumption that the initial covariance tensor is given in the form 
$F_{0}=F(0,x,y,z) = \rho(z) \, \mbox{Id}$.  By assuming the fluid is moving in 
one direction, the anisotropic LANS-$\alpha$ equations reduce
to a one-dimensional partial differential equation for the mean velocity.  
The authors use the Galerkin method to obtain the existence of weak solutions 
to the anisotropic equations.

In this chapter, we use the anisotropic LANS-$\alpha$ equations as a numerical 
model for two classical examples of laminar velocity profiles.   Unlike 
other approaches to the modeling of tubulence, the Lagrangian averaging 
approach allows us to capture the large scale motion of a fluid in laminar 
regimes.  The 
first example that we consider is the Poiseuille flow in the channel and pipe 
domains.  Laminar Poiseuille flow occurs when an incompressible fluid 
with no-slip boundary conditions is driven by a constant upstream pressure 
gradient, yielding a symmetric parabolic stream-wise profile.   We assume 
the velocity $u$ satisfies the steady Navier-Stokes equations and we 
numerically solve
the anisotropic equations for the initial covariance tensor $F_{0}$.  In 
particular, we calculate the matrix $F$ such that pair $(u,F)$ is a 
solution to the 
anisotropic LANS-$\alpha$ equations.  The degeneracy rate of our 
numerical solution $F$ near the wall of the channel is in good agreement with the logarithmic decay rate 
given by Coutand and Shkoller \cite{cs03a}.  In addition, we show that in the 
boundary layer, $||F(t, \cdot)||_{L^{\infty}(\Omega)}$ does not remain 
bounded as $t \rightarrow \infty$, answering, at least numerically, a question 
stated in \cite{cs03a}.  

In the last section of the chapter, we study shear flow solutions to 
the anisotropic LANS-$\alpha$ equations.  Shear flow occurs in the channel 
when one side of the boundary is moving while the other side remains fixed.
The derivation by Marsden and Shkoller \cite{ms03} relied 
upon no-slip boundary conditions for the mean velocity field $u$.  In order to 
find a solution to the inhomogeneous problem, we use the classical 
technique of introducing a new vector field 
that vanishes on the boundary and that solves the anisotropic equations with 
additional forcing terms.  For the domains and the initial data that submit a 
unique solution to the new formulation, we are able to solve for the 
unknown mean velocity.  In Section \ref{inhomog}, we show that shear flow 
velocity solutions to the anisotropic LANS-$\alpha$ equations exist if the 
initial covariance tensor is not required to be positive.  This turns out 
to be a natural assumption and we compute the shear flow solutions directly.

\subsection{The covariance tensor $F(t,x).$}

\begin{figure}[tbp]
    \centering
    \resizebox{10cm}{!}{\includegraphics{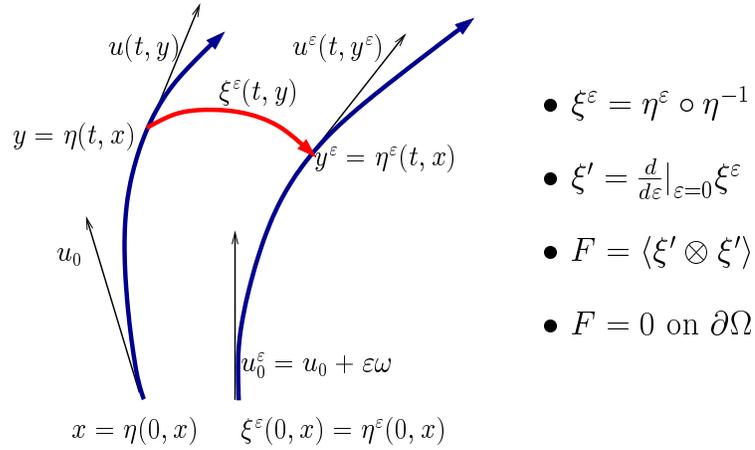}}
    \caption{The construction of $F$.}
    \label{fig:xi}
\end{figure}

The matrix $F$ is defined by Marsden and Shkoller 
\cite{ms03} 
as the ensemble average of the tensor product of the Lagrangian 
fluctuation vector.  Specifically, if $\eta$ and $\eta^{\varepsilon}$ 
are the particle trajectories of the fluid velocity $u$ and averaged 
fluid velocity $u^{\varepsilon}$, respectively, then the Lagrangian 
fluctuation vector 
$\xi^{\varepsilon} := \eta^{\varepsilon} \circ \eta^{-1}$ (see Figure 
\ref{fig:xi}) and 
\[ F := \displaystyle{ \left< 
\frac{d}{d\varepsilon}\Big|_{\varepsilon=0} \xi^{\varepsilon} \otimes 
\frac{d}{d\varepsilon}\Big|_{\varepsilon=0} \xi^{\varepsilon} \right> }\]
where $<\cdot>$ denotes the average over all possible solutions 
$u^{\varepsilon}$ (see \cite{ms03}).  
Because these fluctuations are necessarily zero along the boundary,
\[ F(t,x) = 0 \,\, \mbox{for} \,\, t \geq 0 \]
on the boundary of the domain.  As noted in \cite{cs03a}, it is 
unknown \emph{a priori} whether  
$||F(t, \cdot)||_{L^{\infty}}$ remains bounded for all time.  
In next section, we show numerically that $||F(t, \cdot)||_{L^{\infty}}$ is 
strictly increasing in time in the viscous boundary layer of the channel and 
pipe under a steady mean flow.

\section{Homogeneous boundary conditions.}\label{homog}

Laminar Poiseuille flow occurs when an incompressible fluid in a straight 
channel, or pipe, is driven by a constant upstream pressure gradient, yielding a 
symmetric parabolic stream-wise velocity profile.  In this section, 
we begin with the steady Poiseuille flow associated with the Navier-Stokes 
equations.  We construct an initial covariance tensor such that the 
steady flow solution and the 
covariance tensor solve the anisotropic equations in the channel.  First, 
we reduce the anisotropic 
LANS-$\alpha$ equations to the channel and pipe under the assumption that the 
initial covariance tensor is a multiple of the identity matrix.

\subsection{Given the initial covariance matrix.}  

\noindent {\bf Channel.}   The three-dimensional channel is given by 
$\Omega = \mathbb{R}^2 \times 
[-h,h]$ with coordinates ${\bf x}=(x,y,z)$.  We shall assume that the velocity
 vector is of the form 
\[ {\bf u}(t,x,y,z) = (u(t,z),0,0), \]
and that the initial covariance matrix is written as 
\[ F_{0}=F(0,x,y,z) = \rho(z) \mbox{Id}, \]
such that $\rho(\pm h) = 0$. 

Suppose $\rho(z)$ is given.  As discussed in \cite{cs03a}, when no-slip 
boundary conditions are prescribed for the mean velocity, the 
nonlinear term vanishes and equation (\ref{aniso}) reduces to the 
following system
\begin{equation}\label{aniso_channel}
    \begin{array}{c}
\partial_t (u-\alpha^2(\rho u')')-\nu((\rho u')' - 
\alpha^2 (\rho(\rho u')'')') = - c, \\
u(t,\pm h) = 0, \hspace{.2in} u(0,z)=u_0(z), 
\end{array}
\end{equation}
where $c :=\partial_x p$ is constant since $\partial_y p= \partial_z p 
=0$.  To address the degeneracy at the boundary, Coutand and Shkoller 
define a weighted Sobolev to set as the functional framework.  This 
space is based on the following function:  For $\epsilon < \delta << 1$ 
define $\rho_{CS}(z) \in C^{\infty} [-h,h]$ to be the positive function
\begin{equation} \label{p_choosen}
    \begin{array}{c}
\rho_{CS}(z) = \left\{ 
\begin{array}{cc}
d \sqrt{|\log d|}, & 0 \leq d(z) \leq \epsilon, \\
1, & d(z) \geq \delta, 
\end{array} 
\right.  \\
d \ \, \mbox{= normalized distance function to the boundary}.
\end{array}
\end{equation}
Coutand and Shkoller prove the existence of a unique global 
weak solution to equation (\ref{aniso_channel}) in a weighted 
Sobolev space.  In addition, they show that the covariance tensor $F$ must 
degenerate like $\rho_{CS}$ in the viscous boundary layer.  In the next section, 
we verify that our numerical solution $\rho$ matches $\rho_{CS}$ near 
the walls of the channel.

\noindent {\bf Pipe.}  The three-dimensional 
pipe is given by $\Omega = \mathbb{R}^3$ with 
coordinates ${\bf x} = (x,y, z)$ such that $0\leq \sqrt{x^2+y^2}\leq a$.  
We shall assume that the velocity vector is radially symmetric of the form
\begin{equation}\label{u_pipe}
{\bf u}(t,x,y,z) = (0,0,u(t,r)),
\end{equation}
where $r:=\sqrt{x^2+ y^2}$ and that the initial covariance matrix 
\[ F_{0}=F(0,x,y,z) = \rho(r) \mbox{Id}.\]
We let ${\bf \eta}$ denote the Lagrangian flow of ${\bf u}$ satisfying the 
initial value problem (\ref{eq:eta}). 
Then definition (\ref{u_pipe}) implies that
\[ \displaystyle{ {\bf \eta}(t,{\bf x}) = \left(x,y, z+\int_0^t u(s,r) \, ds 
\right).} \]
Therefore
\begin{equation} \label{F_pipe}
\displaystyle{ F(t,{\bf x}) := D{\bf \eta}(t,{\bf x}) \cdot F_0({\bf x}) \cdot 
D{\bf \eta}(t,{\bf x})^T = \rho(r) \left[
\begin{array}{ccc}
1 & 0 & \frac{x}{r}\mathcal{U} \\
0 & 1 & \frac{y}{r}\mathcal{U} \\
\frac{x}{r}\mathcal{U} & \frac{y}{r}\mathcal{U} & 1+\mathcal{U}^2
\end{array}
\right],} 
\end{equation}
where $\displaystyle{\mathcal{U}(t,{\bf x}):=\int_0^t \partial_r u(s,r)\,ds}$. 
Together with definition (\ref{u_pipe}), $C{\bf u} := \operatorname{div}(\nabla u \cdot 
F)$ reduces to 
\[ C{\bf u} = (0,0,\frac{1}{r}\left(r \rho u'\right)') \]
with $u'(t,r):=\partial_r u(t,r)$.  The nonlinear term in the 
anisotropic LANS-$\alpha$ vanishes and the equations reduce to the following 
system
\begin{equation}\label{aniso_pipe}
    \begin{array}{c}
\displaystyle{\partial_t (u- \alpha^2\frac{1}{r}(r \rho u')')- 
\frac{\nu}{r}\left( (r\rho u')'-\alpha^2  (r\rho 
(\frac{1}{r}(r\rho u')')')'\right) = - c,} \\
u(t,a) = 0, \hspace{.2in} u(0,r)=u_0(r),  
\end{array}
\end{equation}
where $c:=\partial_z p $ is constant since $\partial_x p= \partial_y
p =0$. 

\subsection{Given the mean velocity profile.}\label{sec:num_F}

In this section, we provide the numerical 
results for steady channel and pipe flow assuming no-slip boundary conditions.  
We start with the steady solution $u$ to the Navier-Stokes equations and proceed 
to find a $\rho$, and therefore a $F$, such that $(u,F)$ solve the anisotropic 
LANS-$\alpha$ equations.

\noindent {\bf Channel.}  We begin by assuming the velocity 
field is the steady flow given as 
${\bf u}(t,x,y,z) = (u(z),0,0)$ that satisfies the no-slip boundary condition 
$u(\pm h) = 0$.  For the Navier-Stokes equations, the classical Poiseuille 
flow is given by 
\begin{eqnarray}\label{poiseuille}
u(z) &=& \beta (h^2 - z^2) \\
p(x) &=& p_0 - 2\nu \beta x, \nonumber
\end{eqnarray}
where $\beta, p_0 >0$.  In order for $(u,\rho)$ to 
solve equation (\ref{aniso_channel}) with the velocity $u(z)$ given 
by equation (\ref{poiseuille}), the function $\rho(z)$ must solve the 
ordinary differential equation
\begin{equation}\label{channel_rho}
    \begin{array}{c}
(z\rho)'-\alpha^2 (\rho (z \rho)'')' =1, \\
\rho(\pm h) = 0 . 
\end{array}
\end{equation}
Solving equation (\ref{channel_rho}) numerically leads us to the following 
Proposition.

\begin{figure}[!ht]
     \centering
     \subfigure[The solution $\rho(z)$ of equation (\ref{channel_rho}). ]{
	  \label{fig:rho}
	  \includegraphics[width=.45\textwidth]{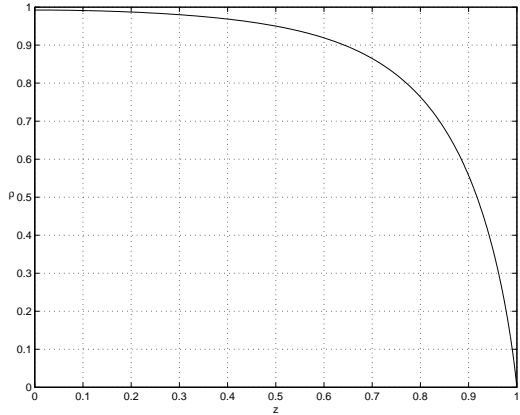}}
    \hspace{.2in}
     \subfigure[The solution $\rho(z)$ and the function $\rho(z)$ used in 
     \cite{ms03} near the boundary.]{
	  \label{fig:CS_rho}
	  \includegraphics[width=.45\textwidth]{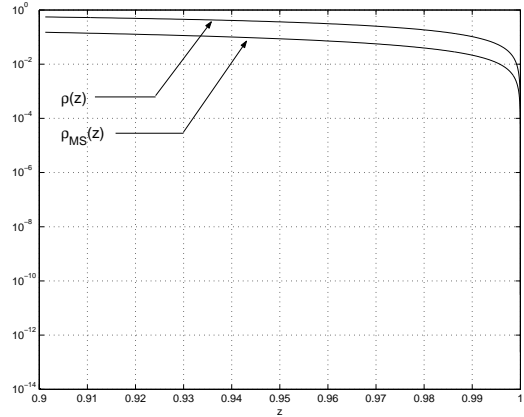}} \\
     \vspace{.2in}
     \subfigure[The eigenvalues of $F(t,z)$ at $t=2$.]{
	   \label{fig:ev_of_F}
	   \includegraphics[width=.45\textwidth]
		{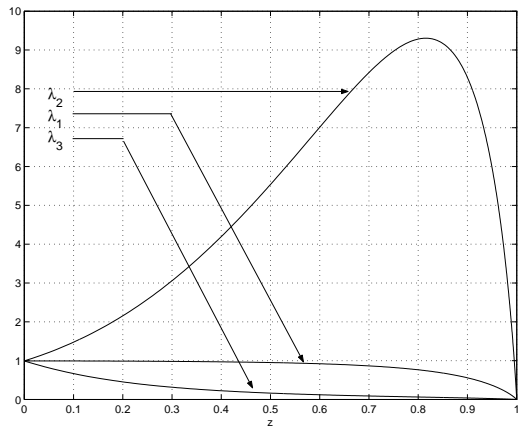}}
    \hspace{.2in}
     \subfigure[The evolution of $||F||_{L^\infty}$.]{
	   \label{fig:F_max}
	  \includegraphics[width=.45\textwidth]{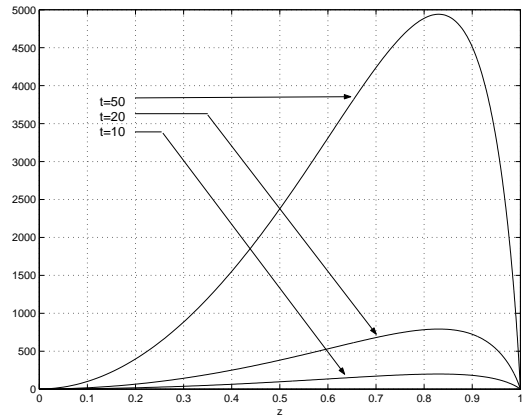}}
     \caption{For $h=\beta =1$ and $\alpha=0.1$. Because of symmetry, 
     the graphs (a), (c), and (d) are on $[0,1]$ only.}
     \label{fig:channel}
\end{figure}

\begin{proposition} Let the function $u$ be the steady solution 
    (\ref{poiseuille}) to the incompressible Navier-Stokes equations in 
    the smooth three-dimensional 
    channel.  The pair $(u,\rho)$ solves the channel anisotropic 
LANS-$\alpha$ equation (\ref{aniso_channel}) when the function $\rho(z)$ 
solves equation (\ref{channel_rho}).  A numerical solution is given in 
Figure \ref{fig:rho} when $\alpha = 0.1$ and $h = \beta = 1$.  In addition, 
the covariance tensor $F$ is defined as 
\begin{equation*}
F(t,z) := \rho (z) D\eta(t,{\bf x}) D\eta(t,{\bf x})^T = \rho(z) \left[
\begin{array}{ccc}
1+4\beta^2 t^2 z^2 & 0 & -2\beta tz \\
0 & 1 & 0 \\
-2 \beta tz & 0 & 1
\end{array}
\right].
\end{equation*}
\end{proposition}

\begin{remark} We assume that $\rho(z)$ is positive and 
    shares the same symmetry across the channel as the velocity $u(z)$.  In 
    addition, we assign in the center of the channel the values 
    $\rho(0) = A$ and $\rho^{\prime}(0)=0$.  By a numerical shooting method, we 
    find a value for $A$ that forces the zero boundary condition $\rho(1)=0$.
    This solution is plotted in Figure \ref{fig:rho}.
\end{remark}

\begin{remark} The numerical solution $\rho$ decays at the same rate as the 
    function $\rho_{CS}$ near the boundary.  In Figure \ref{fig:CS_rho}, 
    we plot both functions near the wall of the channel.
\end{remark}

To illustrate the dynamics of $F$ in the channel, we calculate 
its eigenvalues and conclude the following Proposition.

\begin{proposition} \label{F_prop}
    The eigenvalues of $F(t,z)$ are
\begin{eqnarray*}
\lambda_1(t,z)&=& \rho(z), \\
\lambda_{2,3}(t,z) &=& \frac{1}{2} \rho(z) \left( 2+\mathcal{U}^2 \pm 
|\mathcal{U}|\sqrt{\mathcal{U}^2+4} \right) = \rho(z) \left(1+2\beta^{2} 
t^2z^2 \pm 2\beta t|z| \sqrt{1+\beta^{2} t^2z^2} \right). 
\end{eqnarray*}
The eigenvalues are plotted in Figure \ref{fig:ev_of_F} for the case of 
$\alpha = 0.1$, $\beta =1$, and the function $\rho$ given as the solution 
to equation (\ref{channel_rho}).  Furthermore, $||F(t,\cdot)||_{L^\infty(\Omega)} = 
O(t^2)$ and although our velocity is steady throughout the channel, the 
covariance tensor 
\[ ||F(t,\cdot)||_{L^{\infty}(\Omega)} \rightarrow \infty \hspace{.2in} 
\mbox{as} \hspace{.2in} t \rightarrow \infty \]
near the boundary in the direction of the flow.  This is illustrated
in Figure \ref{fig:F_max}.
\end{proposition}


\newpage

\newpage

\noindent {\bf Proof.} The eigenvalues are straightforward to compute.  Since 
the largest eigenvalue 
$\lambda_2(\cdot, x) = O(t^2)$ and $\lambda_2 \leq ||F(t,z)||_{L^\infty}$, 
then $||F(t, \cdot)||_{L^\infty} = O(t^2)$.  The vector
\begin{equation*}
v(t,z) = C \left[
\begin{array}{c}
tz + \sqrt{1+t^2 z^2} \\
0 \\
-1
\end{array} \right]
\end{equation*}
solves the eigenvalue equation $Fv=\lambda_2v$.  Therefore, as 
$t \rightarrow \infty$,
the norm of $F$ increases to infinity in the direction of the flow 
($x$-direction).
\qed

\noindent {\bf Pipe.}  The format and conclusions of this subsection 
follow closely the results in the channel.  We begin by assuming the velocity 
field is a steady velocity field given as 
${\bf u}(t,{\bf x}) = (0,0,u(r))$ that satisfies the no-slip boundary 
condition $u(\pm a) = 0$.  For the Navier-Stokes equations, the classical 
pipe flow is given by 
\begin{equation}\label{pipe}
    \begin{array}{c}
u(r) = \beta (a^2 - r^2) \\
p(r) = p_0 - 4\nu \beta r, 
\end{array}
\end{equation}
where $\beta, p_0 >0$.  In order for $(u,\rho)$ to 
solve equation (\ref{aniso_pipe}) with the velocity $u(r)$ given 
by equation (\ref{pipe}), the function $\rho(r)$ must solve the 
ordinary differential equation

\begin{equation}\label{eq:pipe_rho}
   \begin{array}{c}
       r(r\rho'+2\rho) - \alpha^2(r\rho (r\rho''+3\rho'))' = 2r \\
       \rho(a) = 0.
\end{array}
\end{equation}
As in the channel domains, we arrive at the following Proposition.

\begin{figure}[thp]
    \centering
    \resizebox{8cm}{!}{\includegraphics{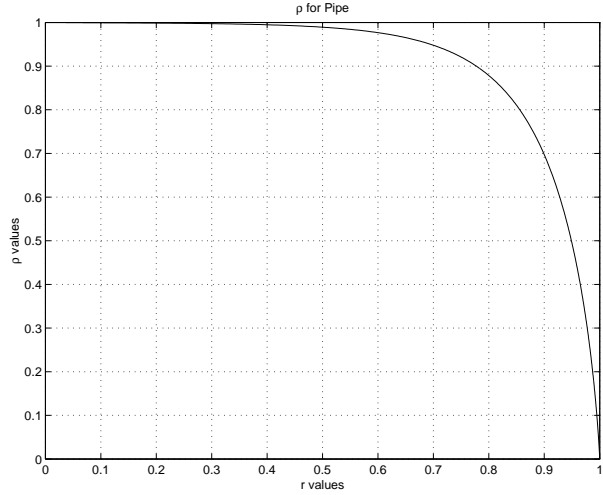}}
    \caption{The solution $\rho(r)$ of equation (\ref{eq:pipe_rho}.)}
    \protect\label{fig:pipe_rho}
\end{figure}

\begin{proposition} Let the function $u$ be the steady solution (\ref{pipe}) 
to the incompressible Navier-Stokes equations in the smooth three-dimensional 
   pipe.  The pair $(u,\rho)$ solves the pipe anisotropic 
LANS-$\alpha$ equation (\ref{aniso_pipe}) when the function $\rho(r)$ solves 
equation (\ref{eq:pipe_rho}).  A numerical solution is given in Figure 
\ref{fig:pipe_rho} when $\alpha = 0.1$, $\beta =1$, and the radius of the pipe is $a = 1$.  
In addition, the covariance tensor $F$ is given by (\ref{F_pipe}) and 
has the same eigenvalues as the channel covariance tensor.  
\end{proposition}

%

%

\subsection{Conclusions.}

As expected, the anisotropy of the fluid is of fundamental importance 
in bounded domains.  We have demonstrated that the dynamics of $F$ are 
just as important in laminar regimes as they are in turbulent regimes.  
Supposing that the covariance tensor is a constant multiple of 
the identity matrix, while accurate for very short time and in the 
center of the channel, becomes an inaccurate assumption in the viscous 
boundary layer.  In both the channel and the pipe geometry, the evolution of 
$F$, illustrated by Figure \ref{fig:ev_of_F}, is non-decreasing in 
time and achieves its greatest value in this boundary layer.   As time 
increases, the location of the maximum value of 
$||F||_{L^{\infty}}$ in the channel approaches the limit $z=0.89$.
Therefore an accurate model for fluid motion in the entire channel or 
pipe should be founded on the anisotropic model, rather than the isotropic 
version, in this boundary region.   

The logarithmic degeneracy rate of the covariance tensor in the channel 
agrees well with the decay rate of the function computed in \cite{cs03a}.  
Since the dynamics of fluid motion in the channel and pipe are similar, it 
is not surprising that the covariance tensor for the channel and pipe are the 
same.  The 
next step in studying the numerical properties of the anisotropic model is 
to solve for the covariance tensor when the velocity is given as a 
time-dependent solution to the Navier-Stokes problem.  This is a necessary 
step in understanding the dynamics of the covariance tensor $F$.

\section{Inhomogeneous boundary conditions.} \label{inhomog}

In many physical models, ranging from the study of blood flow to the 
modeling of earthquakes, at least one of the boundary components of 
the fluid container are in motion.  In these cases, 
the mean velocity of the fluid at the wall will no longer be zero.  In this 
section, we study the anisotropic LANS-$\alpha$ equations with 
inhomogeneous boundary data.  As a specific example, we consider the mean 
fluid motion in a channel when one of the boundary walls is not fixed.  
This motion is called shear flow.  

\subsection{The LANS-$\alpha$ equations with inhomogeneous boundary 
data.}

Let $\Omega$ be a bounded fluid container in $\mathbb{R}^{n}$ with 
boundary $\partial \Omega$ and suppose that the mean velocity field 
$u=g \neq 0$ on $\partial \Omega$.  To find a solution pair $(u,F)$ 
that solves the anisotropic LANS-$\alpha$ equations with
inhomogeneous boundary conditions, we choose a divergence-free vector 
field $v(t, \cdot) \in C^{\infty}(\Omega)$ such that $v=g$ on the 
boundary.  With $v$ given, we define a new vector field $w=u-v$.  The 
vector field $w$ is zero on the boundary and solves the anisotropic 
LANS-$\alpha$ equations with $u$ replaced with $w+v$.  Namely, we now 
search for a solution to the following system of partial differential 
equations
\begin{equation}\label{aniso_inhom}
    \begin{array}{c}
(1-\alpha^2 C) \left( \partial_t w - \nu \mathbb{P} Cw +\operatorname{div} 
(v \otimes w + w \otimes v) \right) = - \operatorname{grad} p - (1-\alpha^2 
C) \left( \partial_t v - \nu \mathbb{P} Cv \right), \\
\operatorname{div} w(t,x) = 0, \hspace{.2in} w=0 \,\, \mbox{on} \,\, 
\partial \Omega,  \\
\partial_t F + \nabla F \cdot w - \left( \nabla w \cdot F + 
[\nabla w \cdot F]^T \right) = \nabla F \cdot v + \nabla v \cdot F + 
[\nabla v \cdot F]^T,  \\
w(0,x) = u_0(x)-v(0,x),  F(0,x)= F_0(x). 
\end{array}
\end{equation}

In the three-dimensional torus $\mathbb{T}^{3}$, the trivial existence 
and uniqueness of a solution to (\ref{aniso_inhom}) follows from the 
smoothness of $v$, Theorem 1 in \cite{ms03}, and Theorem 
\ref{stokes_wp} in Chapter \ref{chapter:spwp}.  This result is stated 
as the following proposition.
\begin{proposition}
    For $s>7/2$, and $w_0 \in 
    H^s_{\operatorname{div}} 
    (\mathbb{T}^3)$, $F_0 \in [H^s_{\operatorname{per}} 
    (\mathbb{T}^3)]^{3 \times 3}$, with $F_0>0$, there exists a unique solution 
    $(w, F)$ with $w \in C^0([0,T]; H^s_{\operatorname{div}}) \cap L^2 (0,T; 
    H^{s+1}_{\operatorname{div}})$ and $F \in C^0([0,T]; 
    [H^s_{\operatorname{per}}]^{3 \times 3})$ to equations 
    (\ref{aniso_inhom}), where $T$ depends on the initial data.
\end{proposition}

In arbitrary bounded domains, it is unknown whether solutions 
exist to the anisotropic equations.  As we demonstrate in the next 
section, under certain limiting conditions we can find a solution to 
equation (\ref{aniso_inhom}).

\subsection{Given the initial covariance matrix.}

To study shear flow velocity solutions to the anisotropic 
LANS-$\alpha$ equations, we need to restrict the full equations to the 
three-dimensional channel by making a number of 
limiting assumptions.  We assume that the initial covariance tensor $F_0 = 
F(0,x,y,z)=\rho(z) \mbox{Id}$ is given and the mean velocity has the 
form ${\bf u} = (u(t,z),0,0)$ for $z\in[-h,h]$, where $u \neq 0$ on the 
boundary.  

As was done in the general case above, we choose a vector 
field $v(t,\cdot) \in C^\infty [-h,h]$ such that $v(t,\pm h) = u(t,\pm 
h)$ for all $t\geq 0$.  The vector field $w=u-v$ then solves the 
following one-dimensional problem
\begin{equation}\label{aniso_channel_inhom}
\begin{array}{c}
    \mathcal{L}^\alpha w = -c +f, \\
w(\pm h) = 0,  w(0,z)=u_0(z)-v(0,z) 
\end{array}
\end{equation}
where $\mathcal{L}^\alpha$ is the linear operator defined as
\begin{equation*}
\mathcal{L}^\alpha u := \partial_t(u-\alpha^2(\rho u')')-\nu((\rho u')' - 
\alpha^2 (\rho (\rho u')'')'), 
\end{equation*}
and $f:= -\mathcal{L}^\alpha v$.  Then ${\bf u}(t,x,y,z) := 
(w(t,z)+v(t,z), 0,0)$ solves equation 
(\ref{aniso_channel}) with inhomogeneous boundary conditions.  

\subsection{Given the mean velocity field.}

Suppose that we are given the steady velocity vector 
field ${\bf u} = (u(z),0,0)$ with $u(z) \neq 0$ on the boundary.  
To find a solution to the anisotropic model, we need to find a 
function $\rho(z)$ such that
\begin{equation*}
    \nu ((\rho u')'-\alpha^2(\rho (\rho u')'')') = -c. 
\end{equation*}
But since the vector field $u$ doesn't satisfy the no-slip boundary 
conditions, we do not know the correct boundary conditions for $\rho$.  
Rather, we choose $v\in C^{\infty}[-h,h]$ with $v=u$ on boundary.  
Then $w=u-v$ and the pair $(w,\rho)$ solve
\begin{equation}\label{channel_inhom}
\mathcal{L}^\alpha w = -c +f, 
\end{equation}
where $f:= - \mathcal{L}^\alpha v$.  Since $w=0$ on the boundary we also have 
the boundary condition $\rho(\pm h) =0$.  
Since $w$ is given we may find 
$\rho$ satisfying equation (\ref{channel_inhom}) with zero boundary 
conditions.  However, unlike the case when $u$ satisfies no-slip boundary 
conditions, the forcing in equation (\ref{channel_inhom}) removes any 
a priori statement about the positivity of $\rho$. 

\noindent {\bf Shear flow solution.}  As an example, we show the existence of 
a shear flow velocity solution to the anisotropic LANS-$\alpha$ equations.  The 
steady shear flow velocity solution 
to the Navier-Stokes equations with the boundary conditions $u(-h)=0$ 
and~$u(h)=U>0$ is 
\begin{equation}\label{eq:u_shear}
    \displaystyle{u(z) = \frac{U}{2h}(z+h)}
\end{equation} 
with the pressure function $p \equiv 0$.  The function 
$$\displaystyle{w(z) = \frac{U}{4h^2}(h-z)(h+z)}$$ is zero on the boundary and 
solves (\ref{channel_inhom}) with $\displaystyle{f = - \mathcal{L}^\alpha 
\left( \frac{U}{4h^2}(z+h)^2 \right)}$ when $\rho(z)$ is a solution of 
\begin{eqnarray*}
\rho' -\alpha^2 (\rho \rho'')' &=& 0, \\
\rho(\pm h)&=&0.
\end{eqnarray*}
Namely,
\[ \rho(z) = \left\{ 
\begin{array}{ll}
\frac{1}{2\alpha^2}(z-a)(z-b) & \mbox{for}\,\, -h \leq a \leq z \leq b \leq h\\
0 & \mbox{otherwise} 
\end{array} \right. . \]
If we want $\rho \in C^1[-h,h]$ then $a=-h$, $b=h$.  We conclude that 
the pair $(u,\rho)$ with $u$ define by equation (\ref{eq:u_shear}) 
and $\displaystyle{\rho(z) = \frac{1}{2\alpha^{2}} (z-h)(z+h)}$ is a shear flow 
solution which solves equation (\ref{aniso_channel}) with inhomogeneous 
boundary conditions.


\newpage
\pagestyle{myheadings}
\markright{   \rm \normalsize BIBLIOGRAPHY.}
\addcontentsline{toc}{chapter}{{\bf Bibliography}}
\bibliographystyle{hplain}
\thispagestyle{myheadings}
\bibliography{jpeirce_dissertation}

\end{document}